\documentclass[11pt,twoside]{article}

\usepackage{amsfonts}
\usepackage{amssymb}
\usepackage{amsmath}
\usepackage{graphicx}
\usepackage{hyperref}
\numberwithin{equation}{section}
\newtheorem{Prop}{\bf Proposition}[section]

\newtheorem{defn}{\bf Definition}[section]
\newtheorem{Rem}{\bf Remark}[section]
\newtheorem{Ex}{\bf Example}[section]
\newtheorem{Th}{Theorem}[section]

\begin{document}
\def \b{\Box}
\def \to{\mapsto}
\def \e{\varepsilon}
\begin{center}
{\Large {\bf The synthetic presentation of the main research directions in groupoid theory }}
\end{center}

\begin{center}
{\bf Gheorghe IVAN}
\end{center}

\setcounter{page}{1}
\pagestyle{myheadings}

{\small {\bf Abstract}. The purpose of this paper is to present a systematic exposition of the  main results obtained in the studies
carried out in groupoid theory.
{}

\section{Introduction}
\indent\indent There are basically two ways of approaching groupoids. The first one is the category theoretical approach. The
second one is algebraically considering them as a particular generalization of the structure of group. Groupoids are like
groups, but with partial multiplication; i.e., only specified pairs of elements can be multiplied and inverses with respect to
the multiplication exist for each element.

The notion of  groupoid was introduced by H. Brandt (\cite{brand}) and it is developed by
P.J. Higgins in \cite{higgi71}. The algebraic structure of groupoid is similar to a group, with the exception that products of
elements cannot always be defined (\cite{west71, brown, mack87, rare01}).

A generalization of Brandt groupoid  has appeared  in the paper of C. Ehresmann (\cite{ehre50}).

In the last four decades, a large number of researchers have paid great attention to the study of groupoids and their applications.

In this context,  team made up of teaching staff from the Faculty of Mathematics and Informatics of the West University of Timi\c soara researched various important topics from the theory of groupoids and its applications within the scientific seminar: {\it Geometry and Topology Seminar} (Seminarul de Geometrie \c si Topologie). The Geometry and Topology Seminar was founded and coordinated by Professor Dr. Dan I. Papuc since 1972. Among the active members of this scientific seminar we mention: Mircea Puta, Dumitru Opri\c s, Gheorghe Ivan, Marian Degeratu, Mihai Ivan and Vasile Popu\c ta. These researchers have investigated various topics related to groupoids, obtaining outstanding results that have been published in specialized journals. More specifically, we list some representative works: M. Puta (\cite{puta5a, puta5b}), D. Opri\c s (\cite{opris5a, opris5b}), G. Ivan (\cite{giv1a, giv5b, giv1b, giv2b}), M. Degeratu (\cite{dege5a, dege1a, mivde18}), M. Ivan(\cite{miv3d, miv2a, miva06, miva13, miva15}) and V. Popu\c ta (\cite{popu4a, popu4b, popu4c}).

Investigating the studies undertaken in the theory of groupoids, the following five research directions related to this important field of mathematics stand out:\\
{\bf 1.} Groupoids (general theory);\\
{\bf 2.} Topological groupoids and Lie groupoids;\\
{\bf 3.}  Finite groupoids;\\
{\bf 4.} Groupoids endowed with additional structures;\\
{\bf 5.} Groupoid applications.

The paper is organized as follows. In Section 2 we present three equivalent definitions of concept of a Brandt groupoid. Also, general notions and results about  Brandt groupoids are discussed. 
Section 3 present: the definitions of Ehresmann groupoid and groupoid morphism, algebraic substructures of a given Ehresmann groupoid, some universal  operations with these groupoids and main results related to them.
 Section 4 is devoted to the topological groupoids and Lie groupoids. The main properties of these two usefull concepts are discussed.
 A brief enumeration of articles discussing concepts and main results concerning finite groupoids, structured groupoids  and applications of groupoids in various fields is given in Section 5.\\[-0.5cm]

\section{ Groupoids: definitions and basic results}

This section  is devoted to giving basic definitions and some results related to Brandt's groupoids.

\markboth{Gheorghe IVAN}{The synthetic presentation of the main research directions
in groupoid theory}

The notion of groupoid from an algebraic point of view was first introduced  by H. Brandt in a 1926 paper \cite{brand}. From this setting a groupoid (in the sense of Brandt)
can be thought of as a generalization of a group in which only certain multiplications are possible and it contains several  units elements. Other equivalent definition
of groupoid and its properties appear in \cite{brown}, where a groupoid is defined as  small category where each morphism is an isomorphism.
It is very important to note that the {\it definition of groupoid} can be presented in the following {\bf three equivalent ways}: from {\it an axiomatic point of view as
 that of group}  (\cite{cdw87, pata18, avmp20}), from {\it an approach in terms of category theory} (\cite{higgi71, mapi22}) and as
 a {\it particular case of  universal algebra}  (\cite{giv1b}).\\

{\bf 2.1.~The axiomatic definition of the concept  of groupoid }

\begin{defn}
 (\cite{mikw88, wein96})~ {\rm A \textit{groupoid $ G $ over} $ G_{0}~$  or a {\it $ G_{0}-$groupoid}  is a pair $ (G, G _0) $ of nonempty sets such that $
G _0\subseteq G  $ endowed with two surjective maps $~\alpha,\beta : G \rightarrow G_{0},~$
a partially binary operation $~\mu :G_{(2)}\rightarrow G, (x,y)~\rightarrowtail~\mu(x,y):=x y,~$ where\\
$ G_{(2)}=\{(x,y)\in G \times G | \beta(x)=\alpha(y)\}~$ and an injective map $~\iota:G \rightarrow G, x \rightarrowtail \iota(x):=x^{-1} $ satisfying the following properties:\\
$(1)~$({\it associativity})~$ (x y) z = x (y z),~$ in the sense that, if one side of the equation is defined so is the other one and then they are equal;\\
$(2)~$({\it identities}) $~(\alpha (x),x),~(x,\beta (x))\in G_{(2)}~$ and $~ \alpha(x) x = x\beta(x) = x;$\\
$(3)~$({\it inverses}) $~(x^{-1},x),~(x,x^{-1})\in G_{(2)}~$ and $~x^{-1} x = \beta(x),~ x x^{-1} = \alpha(x).~$}
\end{defn}
The elements of $~G_{(2)}~$ are called {\it composable pairs} of $~G.~$ The element $~\alpha(x)~$ (resp. $~\beta(x)~$) is the {\it left unit}  (resp. {\it right unit})
of $~x\in G.~$ The subset $~G_{0} = \alpha(G) = \beta(G)~$ of $~G~$  is called the {\it unit set} of $~G.~$

 A groupoid $ G $ over $~G_{0} $ with the \textit{structure functions} $ \alpha $ ({\it source}), $ \beta $ ({\it target}), $ \mu $ ({\it
partial multiplication}) and $\iota$ ({\it inversion}), is denoted by $ (G, \alpha, \beta, \mu, \iota, G _0) $ or $ ~(G, G_{0}).$

A group $~{\cal G}~$ having $~e~$ as unit element is just a groupoid over $~\{ e\},~$ and conversely, every groupoid with one unit element is a group.

\begin{Ex}
{\rm Let $~{\bf R}^{\ast}~$ be the set of nonzero real numbers and $~a, b\in {\bf R}^{\ast}~$ such that $~a b=1. $ Consider the sets  $ G = {\bf R}^{\ast}\times {\bf R}^{\ast}: = {\bf R}^{\ast 2}~ $ and $ G_{0} =\{ (a, a x)\in G~|~(\forall) x\in  {\bf R}^{\ast}\}. $  It is easy to see that $~ y=ax ~$ if and only if $ x =b y $ for all $~(x, y)\in G. $ Then
$ (G, \alpha, \beta,  \odot, \iota, G_{0})~$ is a groupoid, where the set $~G_{(2)}~$ and  its structure functions are given by:\\
$\alpha : G~\rightarrow~ G_{0}, (x,y)~\rightarrowtail \alpha (x,y):=(x,ax),~\beta : G~\rightarrow~ G_{0}, (x,y)~\rightarrowtail \beta (x,y):=(by,y),~$;\\
$ G_{(2)} = \{ ((x,y), (z, u))\in G \times G ~|~ \beta (x,y) =\alpha (z, u)\} = \{ ((x,y), (z, u))\in G \times G ~|~ z=by \};$\\
$(x,y)~\odot~(by,u):=(x, u),~ (\forall)~x, y, u \in {\bf R}^{\ast}~~$ and $~~\iota(x,y):=(by, ax), (\forall)~x, y\in {\bf R}^{\ast}. $

It is easy to verify that the conditions $ (1) - (3) $ of Definition 2.1 hold. First, we consider $~ (x, y), (z, u), (v, w) \in G.~$   The product
$ (x, y)\cdot (z, u) \cdot (v, w) $ is defined if and only if $~z=by~$ and $~v=bu.~$ We have:\\
$ (1)~~~((x, y)\odot (by, u))\oplus (bu, w)=((x, u)\odot (bu, w) = (x,w)= (x, y)\odot ((by, u)\oplus (bu, w)); $\\
$ (2)~~~\alpha (x, y)\odot (x, y))= (x, ax)\odot (x, y)= (x, t)\odot (bt, y)=(x, y) ~$ and $~(x, y)\odot \beta(x,y)=(x,y)\odot (by, y)= (x,y);$\\
$ (3)~~~(x, y)\odot \iota(x, y)= (x, y)\odot (by, ax)= (x, ax)=\alpha (x,y)~$ and $~\iota(x,y)\odot (x,y)= (by, ax)\odot (x,y)==(by, t)\odot (bt,y)=(by,y)=\beta(x, y).$

This groupoid is denoted with $~{\bf R}^{\ast 2}(a,b).~$ Therefore, for each pair $~(a,b)\in {\bf R}^{\ast 2} ~$   satisfying relation $~ ab=1,~$ an example of groupoid is obtained in this way. We thus have an infinity of groupoids associated with $~{\bf R}^{\ast 2}.~$}\\[-0.5cm]
\end{Ex}
\begin{defn}
 {\rm  Let $(H,H_{0})$ be a pair of  non-empty subsets, where  $~H\subset G$ and $ H_{0}\subset G_{0}.~$\\
 $(i)~$  The pair $(H,H_{0})$  is called {\it  subgroupoid} of  groupoid $~(G, G_{0}),~$ if it is closed under multiplication (when it is defined)  and inversion, i.e.
the following conditions hold:\\
$(1)~~~(\forall) ~x,y\in H~$ such that $~x y~$ is defined, we have $~x y\in H;~$\\
$(2)~~(\forall) ~x\in H ~~\Rightarrow ~~x^{-1}\in H.$\\
$(ii)~$ If $~H_{0}= G_{0},~$ then $~(H,G_{0})~$ is called  {\it wide subgroupoid} of $~(G, G_{0}).~$\\
$(iii)~$ By a {\it normal subgroupoid} of a groupoid $~(G,  G_{0}),~$ we mean a wide subgroupoid $~(H, G_{0})~$ of $~(G, G_{0})~$ satisfying the property : $~$ for all $~x\in G~$ and $~h\in H~$ such that the product $~x h x^{-1}~$ is defined, we have $~x h x^{-1}\in H.~$}
\end{defn}
\begin{defn}
~ {\rm $~(i)~$ Let $~( G, \alpha, \beta, \mu, \iota , G_{0})~$ and $~( G^{\prime}, \alpha^{\prime}, \beta^{\prime}, \mu^{\prime}, \iota^{\prime}, G_{0}^{\prime})~$  be two groupoids.
A  {\it morphism between these groupoids} is a pair $~(f, f_{0})~$ of maps, where $~f: G \longrightarrow G^{\prime}~$ and $~f_{0} : G_{0}~\to~G_{0}^{\prime},~$ such that the
following two conditions are satisfied:\\
$(1)~~~ f(\mu(x,y)) = \mu^{\prime}(f(x), f(y)),~~~(\forall) ~(x,y)\in G_{(2)};~$\\
$(2)~~~ \alpha^{\prime}\circ f = f_{0}\circ \alpha ~$ and $~ \beta^{\prime}\circ f = f_{0}\circ \beta.~$

$(ii)~$ A morphism of groupoids $~(f, f_{0})~$ is said to be {\it isomorphism of groupoids}, if $~f~$ and $~f_{0}~$ are bijective maps.}
\end{defn}
\begin{Ex}
~ {\rm $(i)~$ Let $~M~$  be a nonempty set.  By {\it a quasipermutation of the set} $~M~$ we mean an injective map from a subset of $~M~$ into $~M.~$

We denote by $~G = {\cal S}(M)~$ or $~G = Inj(S)~$ the set of all quasipermutations of $~M,~$ i.e. $~{\cal S}(M) = \{ f~|~ f:A\to M, f~ \hbox{is injective and } \emptyset \neq A\subseteq M~\}.~$

For $~f\in {\cal S}(M),~$ let D(f) be the domain of $~f, ~R(f)=f(D(f))~$ and

$~G_{(2)} = \{~(f,g)~|~ R(f)=D(g)~\}.~$  For $~(f,g)\in G_{(2)}~$ we define $~\mu(f,g) = g\circ f.~$

If $~Id_{A}~$ denotes the identity on $~A,~$ then $~G_{0} = \{~ Id_{A}~|~ A\subseteq M~\}~$ is the set of units of $~G,~$ denoted by $~{\cal S}_{0}(M)~$  and $~f^{-1}~$ is the inverse function from  $~R(f)~$ to $~D(f).~$ The maps $~\alpha,~\beta~$ are defined by  $~\alpha(f) = Id_{D(f)},~ \beta(f) = Id_{R(f)}.~$ Thus $~{\cal S}(M)~$ is a groupoid over $~{\cal S}_{0}(M).~~{\cal S}(M)~$ is called the {\it symmetric groupoid of} $~M~$ or the {\it groupoid of quasipermutations of} $~M.~$

$(ii)~$ For a given groupoid $~( G, G_{0}),~$ let $~({\cal S}(G),{\cal S}_{0}(G))~$ be the symmetric groupoid of the set $~G,~$ where $~{\cal S}_{0}(G)=\{ Id_{A}~|~ A\subseteq G \}.~$

We consider now the set $~{\cal L}(G) = \{~L_{a}~|~a \in G~\}~$  of all left translations\\
 $L_{a}: G ~\longrightarrow~ G, ~x \to~ L_{a}(x) = a x,~$ whenever $~(a,x)\in G_{(2)}.~$

We have $~D(L_{a})= \{~x\in G~|~(a,x) \in G_{(2)}~\}\neq \emptyset,~$ since $~(a,\beta (a))\in G_{(2)}~$ and so $~L_{a}\in {\cal S}(G).~$ Hence  ${\cal L}(G)~$ is a subset of $~{\cal S}(G).~$

For all $~a,b,x\in G~$ such that $~\beta(a)=\alpha(b)~$ and $~\beta(b)=\alpha(x)~$ we have\\
 $~L_{a}( L_{b}(x)) = L_{a}(b x) = a(b x) = (a b) x = L_{a b}(x)~$ and we note that $~ L_{a} \circ L_{b} = L_{a b}~$ if $~(a,b) \in G_{(2)}.~$ Consequently, we have $~ L_{\alpha(x)} \circ L_{x} = L_{x} \circ L_{\beta(x)} = L_{x},~(\forall)~x\in G.~$

For all $~u\in G_{0}~$ we have $~L_{u} = Id_{D(L_{u})},~$ hence $~L_{u}\in {\cal S}_{0}(G)~$  and\\
${\cal L}_{0}(G) = \{ L_{u} ~|~ u\in G_{0}\}~$ is a subset of $~{\cal S}_{0}(G).~$ Since $~{\cal L}(G) \subseteq {\cal S}(G)~$ and the conditions of Definition 2.2 are satisfied, it follows that $~{\cal L}(G)~$ is a subgroupoid of $~{\cal S}(G).~$

This groupoid is called the {\it groupoid of left translations of} $~G.~$}
\end{Ex}

The groupoid of left translations of a given groupoid plays an important role in the proof of Cayley's theorem for groupoids (Theorem 2.1).

\begin{Th} (\cite{giv1b})
Every groupoid $~(G, G_{0})~$ is isomorphic to a subgroupoid of the symmetric  groupoid $~({\cal S}(G), {\cal S}_{0}(G)).~$
\end{Th}

{\bf Proof.} Let $~(G, G_{0}) $ be a groupoid and  $~({\cal L}(G) ,{\cal L}_{0}(G))~$  its groupoid of left translations. It is proved that the pair of functions $~(\varphi, \varphi_{0}), $
where $~\varphi: G~ \longrightarrow~ {\cal L}(G),~\varphi(a) = L_{a},~(\forall)~a\in G~$ and $~\varphi_{0}: G_{0} \longrightarrow~ {\cal L}_{0}(G),~\varphi(u) = L_{u},~(\forall)~u\in G_{0}~$,
 is an isomorphism of groupoids. \hfill$\Box$\\

{\bf 2.2. The concept of groupoid as a small category}

A {\bf groupoid} is a small category $~{\cal G}~$ in which every morphism is an isomorphism (i.e. invertible).More exactly,  the elements of set $~{\cal G}~$   are called {\it morphisms} ( or {\it arrows})  and the elements of  set $~{\cal G}_{0}~$  are called {\it objects} (or {\it vertices}). If  $~f : x\rightarrow y~$ is a morphism from $~{\cal G},~$ then $ x $ is called the {\it domain} (or {\it source}) {\it of} $~f,$  denoted with $~d(f),$ and $ y~$ is called the {\it range} (or {\it target}) {\it of} $ f,$ denoted with $ r(f).$
For each pair of objects $~(x,y)\in  {\cal G}_{0}\times {\cal G}_{0},$ denote  by  $~ {\cal G}(x,y), $ the set  (possibly empty) of morphisms  from $ x $ to $ y,$ that is
$ {\cal G}(x,y)=\{ f\in {\cal G}~|~f : x \rightarrow y\}.$ Then, ${\cal G} =\cup_{x,y\in {\cal G}_{0}}{\cal G}(x,y).~$

For each object $~x\in {\cal G}_{0},~$ the set of morphisms $~{\cal G}(x,x)~$ is denoted by $~{\cal G}(x),~$ i.e. $~{\cal G}(x)=\{g \in {\cal G} ~|~ g: x \rightarrow x \}.$ The identity morphism from $~{\cal G}(x)~$  is $~Id_{x}: x\rightarrow x.$ We identify any object $~x $ of $~{\cal G}_{0}~$ with its identity morphism, that is, $~x = Id_{x} $ and thus $~{\cal G}_{0}\subseteq ~{\cal G}.$ The composition of morphisms of a groupoid $~{\cal G}~$ will be denoted via concatenation.

For the pair $~({\cal G}, {\cal G}_{0} )~$  we define the set $~{\cal G}_{(2)}~$ of composable pairs  and the structure functions $~d, r, \circ, \varepsilon, \iota,~$ given by:\\
$- ~~ d, r: {\cal G}\rightarrow {\cal G}_{0}~$ called respectively {\it source} and {\it target} maps, were $~d(f):=x ~$ and $~ r(f):=y ~$ for all $~f\in {\cal G}(x,y);$\\
$- ~~{\cal G}_{(2)}:=\{ (f,g)\in {\cal G}(x,y)\times {\cal G}(z,u)~|~r(f)=d(g)\};$\\
$-~~$ the {\it composition law} of morphisms  $~\circ : {\cal G}_{(2)} \rightarrow {\cal G}, $ given by $~(f,g) \rightarrow g\circ f;$\\
$-~~$ the {\it object map}  $~\varepsilon: {\cal G}_{0} \rightarrow {\cal G}, $ given by $~x \rightarrow \varepsilon(x):=Id_{x};$\\
$-~~$ the {\it inversion map}  $~\iota: {\cal G} \rightarrow {\cal G}, $ given by $~f \rightarrow \iota(f):=f^{-1}, ~(\forall) f\in {\cal G}.$

\begin{defn} {\rm  A {\it groupoid $~{\cal G}~$ over $~{\cal G}_{0} $ } is a pair $~({\cal G}, {\cal G}_{0} )~$ endowed with the structure functions $~d, r, \circ, \varepsilon~$ and $~\iota, ~$ satisfying the following properties:\\
{\bf (1)} $~ (h\circ g)\circ f = h\circ (g\circ f)~ (\forall) f, g, h \in {\cal G}~$ such that $~r(f)=d(g) $ and $~r(g)=d(h);$\\
{\bf (2)} $~f\circ \varepsilon (d(f)) = \varepsilon (r(f))\circ f =f, $ that is  $~f\circ Id_{d(f)} = Id_{r(f)}\circ f = f~$ for all $ f\in {\cal G};$\\
{\bf (3)} for all $ f\in {\cal G}~$  there exists $~ \iota(f)\in {\cal G}~$ such that $~ f\circ \iota(r(f))=\varepsilon(r(f))~$ and $~\iota(f)\circ f= \varepsilon (d(f)), $ that is  $~f\circ f^{-1}=Id_{r(f)}~$ and $~f^{-1}\circ f= Id_{d(f)}.$}
\end{defn}

A groupoid $~{\cal G}~$ over $~{\cal G}_{0} $  is denoted by $~( {\cal G}, d, r, \circ, \varepsilon, \iota, {\cal G}_{0} ) $ or $ ~({\cal G}, {\cal G}_{0}).$

A groupoid where there is only one object is a usual group.

The {\it isotropy group} associated to an object $ u\in {\cal G}_{0} ~$ is the group $~{\cal G}(u)=\{f\in {\cal G}~|~ d(f) = r(f) = u \}.$

\begin{Prop} Let $~({\cal G}, {\cal G}_{0})~$ be a groupoid and $~f,g\in {\cal G}.~$ The following assertions hold.\\
(i) If there exists $~g\circ f, $ then $~d(g\circ f) = d(f)~$ and $~r(g\circ f) = r(g).$\\
(ii) There exists $~ g\circ f)~$ if and only if there exists $~ \iota(f)\circ \iota(g)~$ and  $~ \iota(g\circ f)= \iota(f)\circ \iota(g).$
\end{Prop}

A {\it subgroupoid} of $~({\cal G}, {\cal G}_{0})~$ is a subcategory $~ ({\cal H}, {\cal H}_{0} ) $ of $~{\cal G}, {\cal G}_{0})~$  which is also a groupoid.
We say that the subgroupoid  $~ ({\cal H}, {\cal H}_{0} ) $ is {\it wide}, if  $~{\cal H}_{0}= {\cal G}_{0}.$ \\
Moreover, a wide subgroupoid $~ ({\cal H}, {\cal H}_{0} ) $ is said to be {\it normal subgroupoid}, if $~f\circ {\cal H}(d(f))\circ f^{-1}\subseteq {\cal H}(r(f)) ~$ or
$~f^{-1}\circ {\cal H}(r(f))\circ f\subseteq {\cal H}(d(f))~$ for any $~f\in {\cal G}.$\\
Recall that  $~f\circ {\cal H}(d(f))\circ f^{-1}= \{f\circ h\circ f^{-1}~|~ f\in {\cal G}~\hbox{with}~d(f)=d(h)=r(h)\}.$\\
The set $~Iso({\cal G}):= \cup_{u\in {\cal G}_{0}}{\cal G}(u)~$ is a wide subgroupoid of the groupoid $~({\cal G}, {\cal G}_{0}),~$ called the {\it isotropy subgroupoid}
 of  $~{\cal G}.$\\
\begin{Rem}
{\rm For the detailed presentation of the groupoid concept in the sense of Definition 2.4  and the fundamental properties related to these groupoids,
 as well as their use in different areas of mathematics, we recommend the most relevant references on groupoids (\cite{higgi71,brown, brown06, mack87, mapi22}).}
\end{Rem}

{\bf 2.3. The Brandt  groupoid as universal algebra}

Let $~(\Gamma, \alpha, \beta, \mu ,\iota )~$ be an universal algebra , where $~\Gamma~$ is a nonempty set  endowed  with
the maps $~\alpha~$ ({\it source})  and $~\beta~$ ({\it target}) $~\alpha,\beta : \Gamma \longrightarrow \Gamma~$ such that $~\alpha(\Gamma)=\beta(\Gamma),$
 the {\it multiplication  map} $~\mu :\Gamma_{(2)}\longrightarrow \Gamma, (x,y)~\longrightarrow ~\mu(x,y),~$ where
$~\Gamma_{(2)} =\{~(x,y)\in \Gamma \times \Gamma~|~\beta (x)=\alpha (y)~\}~$ and
the {\it inversion map} $~\iota:\Gamma \longrightarrow \Gamma, x\longrightarrow \iota(x).~$
We write sometimes $~x\cdot y~$ or $~x y~$ for $~\mu(x,y)~$ and $~x^{-1}~$ for $~\iota(x).~$
The elements of $~\Gamma_{(2)}~$ are called {\it composable pairs} of $~\Gamma.~$

\begin{defn}
{\rm  A ({\it Brandt}) {\it groupoid} is an universal algebra $~(\Gamma, \alpha, \beta, \mu, \iota )~$ which
satisfy the following conditions :\\
$(1)~$({\it associativity})~$ (x y) z = x (y z),~$ in the sense that, if one side of the equation is defined so is the other and then they are equal;\\
$(2)~$({\it identities}) $~(\alpha (x),x),~(x,\beta (x))\in\Gamma_{(2)}~$ and $~ \alpha(x) x = x\beta(x) = x;$\\
$(3)~$({\it inverses}) $~(x^{-1},x),~(x,x^{-1})\in\Gamma_{(2)}~$ and $~x^{-1} x = \beta(x),~ x x^{-1} = \alpha(x).~$}
\end{defn}

The element $~\alpha(x)~$ ( resp., $~\beta(x)~$ ) is the {\it left unit} ( resp., {\it right unit} )  of $~x\in \Gamma.~$ The subset $~\Gamma_{0} = \alpha(\Gamma) = \beta(\Gamma)~$
 of $~\Gamma,~$  denoted sometimes by $~U_{\Gamma},~$ is called the {\it unit set} of $\Gamma$ and we say that
$~\Gamma~$ is a $~\Gamma_{0}~$ - {\it groupoid} or that $~\Gamma~$ is a {\it groupoid over } $~\Gamma_{0}.~$

A $~\Gamma_{0}~$-groupoid $~\Gamma~$ will be denoted by $~(\Gamma,\alpha,\beta, \mu , \iota , \Gamma_{0})~$  or $~(\Gamma, \alpha, \beta ,\Gamma_{0})~$ or $~(\Gamma, \Gamma_{0}).~$
The maps $~\alpha, \beta, \mu ~$ and $~\iota~$ of a groupoid $~(\Gamma , \Gamma_{0})~$ are called the {\it structure functions} of $~\Gamma.~$

A $~\Gamma_{0}~$-groupoid $~\Gamma~$ is said to be {\it transitive}, if the map\\
 $~\alpha\times\beta: \Gamma\longrightarrow \Gamma_{0}\times \Gamma_{0} , ~x~\to~ (\alpha\times \beta)(x) = (\alpha(x),\beta(x))~$ is surjective.
The map $~\alpha\times \beta~$ is called the {\it anchor} of $~\Gamma.~$

In the following proposition we summarize some properties of the structure functions of a groupoid obtained directly from definitions.

\begin{Prop} {\rm \cite{giv1b})}
~For any groupoid $~(\Gamma, \alpha, \beta, \mu, \iota ,\Gamma_{0})~$ the following assertions hold:\\
$(i)~~~\alpha(u) = \beta(u) =\iota(u)= u ~$ and $~u\cdot u = u~$ for all $~u\in \Gamma_{0};~$\\
$(ii)~~~\alpha(x y) = \alpha(x)~$ and $~ \beta(x y ) = \beta(y), ~(\forall)~ (x,y)\in \Gamma_{(2)};~$\\
$(iii)~~~ \alpha \circ i = \beta, ~~ \beta\circ i = \alpha ~$ and $~ i \circ i = Id_{\Gamma};~$\\
$(iv)~~~$({\it cancellation law} ) If $~x z_{1} = x z_{2}~$ (resp., $~z_{1} x = z_{2} x $), then $~z_{1} = z_{2};~$\\
$(v)~~~(x,y)\in \Gamma_{(2)}~~\Longrightarrow~~(y^{-1},x^{-1})\in \Gamma_{(2)}~~$ and $~~(x y )^{-1} = y^{-1} x^{-1};~$\\
$(vi)~~~$ For each $~u\in \Gamma_{0},~$ the set  $\Gamma(u)=\alpha^{-1}(u)\cap\beta^{-1}(u)=\{ x\in\Gamma~|~\alpha(x)=\beta(x)=u~\}$ is a group under the restriction
of $\mu~$ to $~\Gamma(u),$ called the {\bf isotropy group  of $~\Gamma $ at} $u;$\\
$(vii)~~~$  For each $~x\in \Gamma, $ the isotropy groups $~\Gamma(\alpha(x))~$ and $~\Gamma(\beta(x))~$ are isomorphic;\\
$(viii)~~$ If $~(\Gamma, \Gamma_{0})~$  is transitive, then the isotropy groups $~\Gamma(u),~u\in\Gamma_{0}$ are isomorphes.
\end{Prop}

 If $~\Gamma~$ is a $~\Gamma_{0}~$- groupoid, then $~Is(\Gamma) = \{ x\in \Gamma ~|~\alpha(x)=\beta(x)~\}= \cup_{u\in \Gamma_{0}}\Gamma(u)~$ is a groupoid, called  the
{\it isotropy group bundle }  of  $~\Gamma.~$     It is easy to prove that $~Is(\Gamma $ is a normal subgroupoid of $~\Gamma,~$ called the  {\it isotropy subgroupoid} of $~\Gamma.~$ \\

\begin{Ex}
{\rm ~(i)~  {\bf  The pair groupoid over a set}. Let $~X~$ be a nonempty set. Then $~ \Gamma = X \times X~$ is a groupoid with respect to rules:  $~\alpha(x,y) = (x,x),~ \beta(x,y) = (y,y),~$ the elements $~ (x,y)~$ and $~(y^{\prime},z)~$
are composable in $~\Gamma~\Leftrightarrow~y^{\prime} = y~$ and we take $~(x,y)\cdot (y,z) = (x,z)~$ and the inverse of $~(x,y)~$ is defined by $~(x,y)^{-1} = (y,x).~$\\
This groupoid  is called the {\it pair groupoid} associated to $~X~$ and it is denoted with $~{\cal PG}(X).~$ The unit space of  $~{\cal PG}(X)~$  is the diagonal $~\Delta_{X} = \{ (x,x)~|~ x\in X\}.~$  The isotropy group $~\Gamma(u)~$ at $~u = (x,x)~$ is the null group.

(ii)~  {\bf  The disjoint union of a family of groupoids}. If $~\{~\Gamma_{i}~|~i\in I~\}~$ is a disjoint family of groupoids, let $~\Gamma~=~\cup_{i\in I} \Gamma_{i}~$  and $~\Gamma_{(2)} = \cup_{i\in I}\Gamma_{i,(2)}.~$  Here, two elements  $~x,y\in \Gamma~$ may be composed if and only if 
 they lie in the same groupoid $~\Gamma_{i}.~$ This groupoid is called the {\it disjoint union of  groupoids} $~\Gamma_{i}, i\in I,~$ and it is denoted by $~\coprod_{i\in I} \Gamma_{i}.~$ The unit set of this groupoid is $~\Gamma_{0} = \cup_{i\in I} \Gamma_{i,0},~$ where $~\Gamma_{i,0}~$ is the unit set of $~\Gamma_{i}.~$

In particular, the disjoint union of the groups $~G_{i}, i\in I,~$ is a groupoid, i.e. $~G = \coprod_{i\in I}{G}_{i},~$ which be called the {\it groupoid associated to family of groups} $~{G}_{i}, i\in I.~$ \\

For example, let $~GL(n;{\bf R})~$ the general linear group of order $~n~$ over $~{\bf R}.~$  We consider the family of groups $~\{~ GL(k;{\bf R})~|~ 1\leq k\leq n~\}.~$
Then the disjoint union $~\bigcup\limits_{k=1}^{n} GL(k;{\bf R})~$ is a groupoid, called the {\it general linear groupoid of order }$~n~$ {\it over} $~{\bf R}~$ and denoted by $~{\cal GL}(n;{\bf R}).~$

Let the map $~f : {\cal GL}(n;{\bf R})~\to~ {\bf R}^{*}~$ defined by $~f(A) = det (A),~$ for all matrix $~A\in {\cal GL}(n;{\bf R}).~$  We have that $~f~$ is a surjective groupoid morphism from the groupoid $~{\cal GL}(n;{\bf R})~$ onto the multiplicative group of reel numbers different from zero.} 
\end{Ex}

\begin{Rem}
~ {\rm  The Definition 2.5 of the Brandt groupoid is essentially the same as the one given by  Coste, Dazord and Weinstein in \cite{cdw87}. There are
various definitions for Brandt groupoids, see \cite{brown06, avmp20}. 
 The basic concepts from groupoid theory  has been defined and studied in a series of papers, for instance:
normal subgroupoids and quotient groupoid  (\cite{pata18, avmp20}), morphisms of Brandt groupoids (\cite{avmp20});
 presentation of some ways of building up new groupoids from old ones, for instance: direct product of a familly of  groupoids (\cite{giv1b}),
  semidirect product  and solvable groupoid (\cite{mapi22}).}
\end{Rem}

\section{ Ehresmann groupoids}

The Brandt groupoids were generalized by C. Ehresmann in \cite{ehre50}. C. Ehresmann added further structures (topological and differentiable as well as algebraic)
to groupoids, thereby
introducing them as a tool in topology and differential geometry. After the introduction of topological and differentiable  groupoids by Ehresmann  in the 1950's, they
have been studied by many mathematicians with different approaches (\cite{ehre50, brown06, giv1a}).

This section deals with the groupoids in the sense of Ehresmann. For more details about Ehresmann groupoids, we refer
the reader to  (\cite{avmp20, giv1b, mack87}}.\\

{\bf 3.1. The Ehresmann  groupoid: definition and results}

\begin{defn}(\cite{wein96}) {\rm   A {\it groupoid $ \Gamma~$ over $ M $} in the sense of Ehresmann  ({\bf Ehresmann groupoid}) is a pair
$( \Gamma, M )$ of nonempty sets endowed with the surjective maps $~\alpha :\Gamma \rightarrow~ M, ~x \rightarrow ~\alpha(x)~$  and $~\beta : \Gamma~\rightarrow~M,~x\rightarrow \beta(x),$  a partially binary operation $~ \mu:\Gamma_{(2)}\rightarrow \Gamma ,~(x,y)\longmapsto  \mu ( x,y):=x\cdot y,~$  where  $~\Gamma_{(2)}:= \{ (x,y)\in \Gamma\times \Gamma | \beta(x) = \alpha(y) \},~$ an injective map $~\varepsilon :M \to \Gamma,\,\,u \to \varepsilon (u),$  and a map $~\iota:\Gamma \to \Gamma,~x \to \iota (x): = x^{ - 1},~$
satisfying the following properties:\\
{\bf (EG1)}$~$ ({\it associativity})$~~(x\cdot y)\cdot z=x\cdot (y\cdot z)~$  for all $~x,y, z\in \Gamma~$ in the sense that if
one of two products $ (x\cdot y)\cdot z $ and $ x\cdot (y\cdot z) $ is defined, then the other product is also defined and they are equals;\\
{\bf (EG2)} $~$ ({\it units}): for each $ x\in \Gamma~$ follows $~(\varepsilon(\alpha(x)),x),~(x,\varepsilon(\beta(x))) \in \Gamma_{(2)}~$ and  we have $~ \varepsilon(\alpha(x)) \cdot x = x \cdot \varepsilon(\beta(x)) = x;$\\
{\bf (EG3)}$~$ ({\it inverses}): for all $~x \in \Gamma~$ follows $~{x,x^{ - 1}}.~({x^{-1},x} \in \Gamma_{{2}} $ and the following equalities hold $~x \cdot x^{- 1} = \varepsilon(\alpha (x))~$ and $~x^{-1} \cdot x = \varepsilon(\beta(x)).$}
\end{defn}

If $~\Gamma~$ is a Ehresmann groupoid over $~M, $  we will sometimes write $~x y~$ for $ \mu(x,y). $ The set $ \Gamma_{(2)} $ is called  the {\it set of composable pairs} of $ \Gamma. $ The element $~\varepsilon(\alpha(x))$ (resp. $~\varepsilon(\beta(x) $)) is the {\it left unit} (resp. {\it right unit}) of $~x\in \Gamma.~$  Also, the sets $~\varepsilon (M)\subseteq \Gamma~$ and $~M = \alpha(\Gamma) = \beta(\Gamma)~$  are called the {\it unit space} and the {\it base}  of  $~\Gamma,~$  respectively.

 A Ehresmann groupoid $ \Gamma $ over $~M~$  with the \textit{structure functions} $ \alpha $ ({\it source map}), $ \beta $ ({\it target map}),
$ \mu $ ({\it multiplication law}), $~\varepsilon $ ({\it inclusion map}) and $ \iota $ ({\it inversion map}), is denoted by $ (\Gamma, \alpha, \beta,  \mu, \varepsilon, \iota, M)~$  or $~(\Gamma, M)~$ or $~\Gamma.~$ Whenever we write a product in a given a groupoid $~(\Gamma, M),~$  we are assuming that it is defined.

A groupoid $~(\Gamma, M)~$ is said to be {\it transitive} (or {\it connected}), if the map
$~\alpha\times\beta: \Gamma\longrightarrow M\times M, ~x~\to~ (\alpha\times \beta)(x) = (\alpha(x),\beta(x))~$ (called  {\it anchor} of $~\Gamma) $ is surjective .

\begin{Rem}
{\rm Let $ (\Gamma, \alpha, \beta,  \mu, \varepsilon, \iota, M)~$ be a Ehresmann groupoid. If $~ M \subseteq \Gamma~$
and $~\varepsilon = Id_{M},~$ then $(G,\alpha ,\beta, \mu ,\iota ,G_{0} )~$ is a Brandt groupoid.}
\end{Rem}
The basic properties of the Ehresmann groupoids are given in the following two propositions.
\begin{Prop}(\cite{giv1b}) If $~(\Gamma, \alpha, \beta,  \mu, \varepsilon, \iota, M)~$ is a Ehresmann groupoid, then:\\
$(i)~$ If $~(x,y)\in \Gamma_{(2)},~$ then $~\alpha (x\cdot y) = \alpha (x)~$ and $~\beta (x\cdot y) = \beta (y);$\\
$(ii)~$  If $~x\in \Gamma,~ $ then $~\beta (x^{-1})=\alpha (x)~$ and $~\alpha (x^{-1}) =\beta (x);$\\
$(iii)~$  If $~(x,y)\in \Gamma_{(2)},$ then $~(y^{-1},x^{-1})\in \Gamma_{(2)}~ $ and $~(x\cdot y)^{-1}=y^{-1}\cdot x^{-1};$\\
$(iv)~~ x^{-1})^{-1}=x, ~~(\forall) x\in\Gamma;$\\
$(v)~~\alpha (\varepsilon ( u ))= u ~$ and $~\beta (\varepsilon ( u )) = u,~~(\forall )~ u \in M;$\\
$(vi)~~\varepsilon (u)\cdot \varepsilon (u)=\varepsilon (u)~$ and $~(\varepsilon (u))^{-1}=\varepsilon (u)~$ for any $~u\in M;$\\
$(vii)~$ For each $ u \in M, $ the set $ \Gamma(u):=\alpha^{-1}(u)\cap \beta^{-1}(u)= \{x\in \in \Gamma | \alpha(x)=\beta(x)= u \}~ $ is
a group, called the {\bf isotropy group of $ \Gamma $ at $ u.$}
\end{Prop}
\begin{Prop} The structure functions of a Ehresmann groupoid  $~(\Gamma, M)~$ have the following properties:\\
$(1)~~~ \alpha \circ \iota = \beta, ~~ \beta\circ \iota = \alpha ~$ and $~\iota \circ \iota = Id_{\Gamma};~$\\
$~(2)~~~ \alpha \circ \varepsilon= \beta \circ \varepsilon= Id_{M}~$ and $~\iota\circ \varepsilon =\varepsilon.$
\end{Prop}
\begin{Prop} Let $~(\Gamma, M)~$ be a Ehresmann groupoid. Then:\\
$(1)~~$ If $~ \alpha (x) = u ~$ and $~\beta (x) = v,~$ then the map $~\omega:\Gamma(u) \rightarrow \Gamma(v),~ z \mapsto \omega(z):= x^{-1}\cdot  z\cdot x ~$ is
a group isomorphism, i.e. $~\Gamma(\alpha(x))~$ and $~\Gamma(\beta(x))~$ are isomorphes;\\
$(2)~~$ If $~(\Gamma, M)~$ is transitive, then the isotropy groups $~\Gamma(u),~u\in M~$ are isomorphes.
\end{Prop}
Let us  we present some substructures in Ehresmann groupoids and several properties of them
which generalize well-known results in group theory.
\begin{defn}
{\rm  Let $~(\Gamma, M)~$ be a Ehresmann groupoid.  A pair of nonempty subsets $ (H, M^{\prime})$
where $ H\subseteq \Gamma $ and $~M^{\prime}\subseteq M,~$ is called {\it subgroupoid} of $~(\Gamma, M),~$ if:\\
{\bf (ESG1)}$~ \alpha(H)=\beta(H)\subseteq  M^{\prime}~$ and $~\varepsilon(M^{\prime})\subseteq H;$\\
{\bf (ESG2)}$~$ for all $~x,y\in H~$ such that $~(x, y)\in \Gamma_{(2)}~~\Longrightarrow ~x y\in H;~$\\
{\bf (ESG3)}$~$ for all $~x\in H~~\Longrightarrow~~x^{-1}\in H.$}
\end{defn}
\begin{defn}
{\rm Let $ (H,M^{\prime})$ be a subgroupoid of the Ehresmann groupoid $ (\Gamma, M).$\\
{\it(i)} $ (H,M^{\prime})$ is said to be a {\it wide subgroupoid} of $ (\Gamma, M),~$ if $~M^{\prime}= M.$\\
{\it(ii)} A wide subgrupoid $ (H, M)$  of $ (\Gamma, M),$  is called {\it normal subgroupoid,} if for all $ g\in \Gamma $ and for all
 $ h \in H $ such that the product $ g\cdot h\cdot {g^ {-1}} $ is defined, we have $ g\cdot h\cdot {g^ {-1}}\in H .$}
\end{defn}
\begin{Ex}
 {\rm A nonempty set $ X $ may be regarded to be a Ehresmann groupoid over $ X, $ called  {\it null groupoid}.
 For this, we take the following groupoid structure:\\
 $\alpha = \beta = \varepsilon= \iota= Id_{X};$\\
 $~x , y\in X~$ are composable $~~\Longleftrightarrow~~ x = y~$ and $~ x\cdot x = x, ~(\forall)  x\in X.$}
\end{Ex}
\begin{Prop}
{\rm(\cite{giv1a})} $~(i)~$ If $ ~(\Gamma, M) ~$  is a Ehresmann groupoid, then $~(\varepsilon(M), M)~$ is a normal  subgroupoid of $~(G,M), $ called {\bf null subgroupoid}) of $\Gamma.$\\
$(ii)~~$  The set $~Is(\Gamma) = \cup_{u\in M}\Gamma(u)\subset \Gamma~$ is a normal subgroupoid of $~\Gamma.~$  $~(Is(\Gamma), M)~$ is called the {\bf isotropy subgroupoid} of $~\Gamma~$  or {\bf isotropy group bundle} associated to  $~\Gamma.$
\end{Prop}
Let us we present three  general constructions of Ehresmann groupoids.

${\bf \bullet}~~~$ {\bf The direct product of two Ehresmann groupoids.}

Let  $~(\Gamma_{j}, \alpha_{j}, \beta_{j}, \varepsilon_{j}, \mu_{j}, \iota_{j}, M_{j})~$  for $~j = 1, 2~$ be two Ehresmann  groupoids. Consider
$~\Gamma = \Gamma_{1}\times \Gamma_{2}~$ the direct product of the sets $~\Gamma_{j},~j=1,2.~$ We give to $~\Gamma~$ a structure of  groupoid as follows.
 The elements $~x = (x_{1}, x_{2})~$ and $~y = (y_{1}, y_{2})~$ from $~\Gamma~$ are composable if and only if $~(x_{j},y_{j})\in \Gamma_{j,(2)}~$ for $~j=1,2 $ and
we take $~(x_{1}, x_{2} )\cdot (y_{1}, y_{2})= ( x_{1} y_{1}, x_{2} y_{2} ).~$ It is easy to verify that $~(\Gamma_{1}\times \Gamma_{2},\alpha_{1}\times\alpha_{2},\beta_{1} \times \beta_{2}, \varepsilon_{1}\times \varepsilon_{2}, \mu_{1}\times \mu_{2}, \iota_{1}\times \iota_{2}, M_{1}\times M_{2})~$ is a Ehresmann groupoid, called the {\it direct product of Ehresmann groupoids} $~(\Gamma_{1}, M_{1})~$ and  $~(\Gamma_{2}, M_{2}).~$

The direct product of two transitive Ehresmann groupoids is also a transitive Ehresmann groupoid.

${\bf \bullet}~~~$ {\bf The Whitney sum of two Ehresmann groupoids over the same base.}

Let $~(\Gamma, \alpha, \beta, \varepsilon, \mu, \iota, M)~$ and $~( \Gamma^{\prime}, \alpha^{\prime}, \beta^{\prime}, \varepsilon^{\prime}, \mu^{\prime}, \iota^{\prime}, M)~$ be two Ehresman groupoids over $~M.$. The set
$~~( \Gamma\oplus \Gamma^{\prime} = \{~(x,x^{\prime})\in \Gamma\times \Gamma^{\prime}~|~ \alpha(x) = \alpha^{\prime}(x^{\prime}), \beta(x) = \beta^{\prime}(x^{\prime})~\}~$
has a natural structure of Ehresmann groupoid over $~M~$ with the following structure functions:\\
$-{\alpha}_{\oplus},{\beta}_{\oplus}:\Gamma \oplus \Gamma^{\prime}\rightarrow M~$ given by $~~\alpha_{\oplus}(x, x^{\prime}): = \alpha (x),~~~\beta_{\oplus}(x, x^{\prime}): = \beta (x);$\\
$-\varepsilon_{\oplus}: M \rightarrow \Gamma \oplus \Gamma^{\prime}~ $ given by $~\varepsilon_{\oplus}(u): =
(\varepsilon (u),\varepsilon^{\prime}(u)),~ (\forall) u \in M;$\\
$-~$ the elements $~(x,x^{\prime})~$ and $~(y,y^{\prime})~$ are composable in $~\Gamma\oplus \Gamma^{\prime}~$ if and only if  $~(x,y)\in \Gamma_{(2)}~$ and $~(x^{\prime}, y^{\prime})\in \Gamma_{(2)}^{\prime}~$ and $~\mu_{\oplus}((x,x^{\prime}), (y,y^{\prime})) = (\mu(x,y), \mu^{\prime}(x^{\prime}, y^{\prime}));~$\\
$~\iota_{\oplus}: \Gamma \oplus \Gamma^{\prime} \rightarrow \Gamma \oplus \Gamma^{\prime}~ $ given by $~\iota_{\oplus}(x, x^{\prime}): = ({\iota}(x), \iota^{\prime}(x^{\prime})).$\\
Then $~(\Gamma\oplus \Gamma^{\prime}, \alpha_{\oplus}, \beta_{\oplus}), \varepsilon_{\oplus}, \mu_{\oplus},\iota_{\oplus},M)~$ is
 a Ehresmann groupoid, called the {\it Whitney sum}  of $~(\Gamma, M)~$ and $~(\Gamma^{\prime}, M).~$\\
The Whitney sum of two transitive Ehresmann groupoids over $~M~$ is also a transitive Ehresmann groupoid over $~M.$

${\bf \bullet}~~~$ {\bf The induced groupoid of a Ehresmann groupoid}

Let $~(\Gamma, \alpha, \beta, \varepsilon, \mu, \iota, M)~$ be a Ehresmann groupoid over $~M~$ and a given map $~f: X \rightarrow M.~$   The set
$~f^{\ast}(\Gamma)=\{(x,y,a)\in X\times X\times \Gamma~|~f(x)=\alpha(a),f(y)=\beta (a)\}~$ has a canonical  structure of Ehresmann groupoid over $~X~$ with the following structure functions:\\
$-\alpha^{\ast},~\beta^{\ast}: f^{\ast}(\Gamma)\rightarrow X, $ given by $~\alpha^{\ast}(x,y,a): =x,~~~\beta^{\ast}(x,y,a): = y;$\\
$-\varepsilon^{\ast}: X \rightarrow f^{\ast}(\Gamma)~ $ given by $~\varepsilon^{\ast}(x): = (x, x, \varepsilon(f(x))),~ (\forall) x \in X;$\\
$-~$ the elements $~(x, y, a)~$ and $~(y^{\prime}, z, b)~$ are composable in $~ f^{\ast}(\Gamma)~$ if and only if  $~y = y^{\prime}~$ and $~(a,b)\in \Gamma_{(2)}~$ and $~\mu^{\ast}((x, y, a ), (y, z, b)) = (x, z, \mu(a,b));~$\\
$~\iota^{\ast}: f^{\ast}(\Gamma) \rightarrow f^{\ast}(\Gamma)~ $ given by $~\iota^{\ast}(x, y,a): = ( y,x, \iota(a)).$

Then $~(f^{\star}(\Gamma),\alpha^{\star},\beta^{\star},\varepsilon^{\star},\mu^{\star}, \iota^{\star}, X)~$ is a Ehresmann groupoid over $~X, $  called the {\it induced groupoid of $~(\Gamma, M)~$ via  $~f: X \to M. $}

The induced groupoid  $~f^{\star}(\Gamma)~$  of a  transitive Ehresmann groupoid $~\Gamma~$ under $~f~$ is also a transitive Ehresmann groupoid.\\

{\bf 3.2. Morphisms of Ehresmann  groupoids}

Let us we present the notion of {\it Ehresmann  groupoid morphism} and prove several properties of them which generalize well-known results in group theory.
\begin{defn}
{\rm (\cite{mack87})  {\it (i)~} Let $~(\Gamma, \alpha, \beta, \varepsilon, \mu, \iota, M)~$ and $~( \Gamma^{\prime}, \alpha^{\prime}, \beta^{\prime}, \varepsilon^{\prime}, \mu^{\prime}, \iota^{\prime}, M^{\prime})~$ be two Ehresmann groupoids. A {\it morphism of Ehresmann groupoids} or {\it Ehresmann  groupoid
morphism} from $ (\Gamma, M) $ into $ (\Gamma^{\prime}, M^{\prime}) $ is a pair $(f, f_{0}) $ of maps, where $ f:\Gamma \to \Gamma^{\prime} $ and $ f_{0}: M \to
M^{\prime}, $ such that the following two conditions hold:\\
$(1)~~~f(\mu(x,y)) = \mu^{\prime}(f(x),f(y))~~~\forall (x,y)\in \Gamma_{(2)};$\\
$(2)~~~\alpha^{\prime} \circ f = f_{0} \circ \alpha,~~~\beta^{\prime}\circ f = f_{0} \circ \beta.$

{\it (ii)~} A groupoid morphism $~(f, f_{0}) : (\Gamma, M)~\to~(\Gamma^{\prime}, M^{\prime})~$ is said to be {\it isomorphism of Ehresmann groupoids}, if $~f~$ and $~f_{0}~$ are bijective maps.}
\end{defn}

If $~M = M^{\prime}~$ and $~f_{0}= Id_{M}~$ we say that $~f : \Gamma ~\to~\Gamma^{\prime}~$ is a $~M~$-{\it morphism of Ehresmann groupoids}.
\begin{Prop}
If $~(f,f_{0}): (\Gamma, M ) \to ( \Gamma^{\prime}, M^{\prime})~$ is a Ehresmann groupoid morphism then:\\
$~~~~~~~~~f \circ \varepsilon = \varepsilon^{\prime}\circ f_{0}~~~$ and $~~~f \circ \iota = \iota^{\prime}\circ f.$
\end{Prop}
\begin{Ex}
{\rm $~(i)~$ If ($~(\Gamma, \alpha, \beta, \varepsilon, \mu, \iota, M)~$ is a Ehresmann groupoid, then $~(\varepsilon(M), M)~$ is a normal subgroupoid of $\Gamma,~$ called  {\it nul subgroupoid}.\\
$(ii)~$ The {\it kernel} of a Ehresmann groupoid morphism $~(f, f_{0}): (\Gamma, M) \rightarrow (\Gamma^{\prime}, M^{\prime})~$ defined by\\
$~~~~~~~~~~~~~~~~$ Ker$(f) =\{ x\in \Gamma~|~ f(x)\in \varepsilon^{\prime}(M^{\prime})\}$\\
is a normal subgroupoid of $~(\Gamma, M ).$}
\end{Ex}
\begin{Prop}
Let $ (f , f_{0}):(\Gamma,\Gamma_{0})\longrightarrow~ (\Gamma ^{\prime}, \Gamma_{0}^{\prime})$ be a groupoid  morphism. Then:\\
$(i)~~$ If $~(H^{\prime}, H_{0}^{\prime})~$ is a subgroupoid of $ (\Gamma^{\prime}, \Gamma_{0}^{\prime}), $  then $~(f^{-1}(H^{^{\prime }}), f_{0}^{-1}(H_{0}^{\prime}))~$ is a
 subgroupoid of $~(\Gamma, \Gamma_{0}).$\\
$(ii)~~$ If $~(N^{\prime}, \Gamma_{0}^{\prime})~$ is a normal  subgroupoid of $ (\Gamma^{\prime}, \Gamma_{0}^{\prime}), $  then $~(f^{-1}(N^{^{\prime }}), \Gamma_{0})~$ is a
normal subgroupoid of $~(\Gamma, \Gamma_{0}) $  such that $~ Ker(f)\subseteq f^{-1}(N^{\prime}).$
\end{Prop}

 In the paper \cite{giv1a} is introduce the concept of  {\it strong  groupoid morphism}. This   special type of  groupoid morphism plays an important role in groupoids  theory, for instance: see the correspondence theorems for Ehresmann groupoids (\cite{giv1a}).

\begin{defn}
{\rm  A {\it  strong groupoid morphism} is a morphism of Ehresmann groupoids
$(f, f_{0}): (\Gamma, M) \rightarrow (\Gamma^{\prime}, M^{\prime})~ $ satisfying the following condition:\\
{\it $~~~~~~~(\forall) x,y \in \Gamma~$  such that $ f(x)\cdot f(y) $ is defined in $ \Gamma^{\prime}~ \Rightarrow~ x\cdot y~$ is defined in $ \Gamma.$}}
\end{defn}
\begin{Ex}
{\rm $~(i)~$ If ($~(\Gamma, \alpha, \beta, \varepsilon, \mu, \iota, M)~$ is a Ehresmann groupoid, then the anchor map $~\alpha\times \beta: \Gamma \to M\times M~$ is a strong morphism of the groupoid $~(\Gamma, M)~$ into the pair groupoid $~{\cal PG}(M)= (M\times M, M).$ 

$~(ii)~$ If $~(f^{\star}(\Gamma), X)~$ is the induced groupoid of $~(\Gamma, M)~$ via $~f: X \to M, $ then $~(f_{\Gamma}, f): (f^{\star}(\Gamma), X) \to (\Gamma, M),~$ where $~f_{\Gamma} (x,y,a):=f(a),~$ is a groupoid morphism, called the {\it canonical morphism of induced groupoid}. The canonical morphism of induced groupoid $~f^{\star}(\Gamma)~$  is not a strong morphism of groupoids (\cite{giv1a})}.
\end{Ex}
\begin{Prop} (\cite{giv1a})
Let $ (f , f_{0}):(\Gamma, \Gamma_{0})\longrightarrow~ (\Gamma^{\prime} , \Gamma^{\prime}_{0})$  be a strong groupoid  morphism. Then:\\
$(i)~~$ If $~(H, H_{0})~$ is an subgroupoid of $ (\Gamma, \Gamma_{0}), $  then $~(f(H), f_{0}(H_{0}))~$ is a
subgroupoid of $ (\Gamma^{\prime}, \Gamma_{0}^{\prime}).$ In particular, $~(Im(f)=f(\Gamma, f_{0}(\Gamma_{0})) ~$ is a subroupoid of  $ (\Gamma^{\prime}, \Gamma_{0}^{\prime}).$ \\
$(ii)~~$ If the function $ f $ is surjective and  $~(N, \Gamma_{0})~$ is a normal subgroupoid of $ (\Gamma, \Gamma_{0}), $  then $~(f(N), \Gamma_{0}^{\prime})~$ is a normal
subgroupoid of $ (\Gamma^{\prime}, \Gamma_{0}^{\prime}).$
\end{Prop}

Let $~(f, f_{0}) : (\Gamma, \Gamma_{0}) \rightarrow (\Gamma^{\prime}, \Gamma_{0}^{\prime}) $  be a groupoid morphism.\\
$(i)~~$ Denote by  $~{\cal \widetilde S}(\Gamma)~$ (resp. $ {\cal \widetilde N}(\Gamma) $) set of subgroupoids (resp., set of normal subgroupoids) of $ \Gamma $
which contains the kernel of $ f,~$ i.e.:\\[-0.3cm]
\begin{equation}
{\cal \widetilde S}(\Gamma):= \{ H~| H~ \hbox{is a {\it subgroupoid} of}~ \Gamma~ \hbox{and}~ ker (f) \subseteq H\},\label{(3.1)}\\[-0.1cm]
\end{equation}
\begin{equation}
{\cal\widetilde N}(\Gamma):= \{ N ~|~ N~ \hbox{is a {\it normal subgroupoid} of}~ \Gamma~ \hbox{and}~ ker (f) \subseteq N\}.\label{(3.2)}
\end{equation}
$(ii)~~$ Denote by  $~{\cal {S}}(\Gamma)~$ (resp. $~{\cal {N}}(\Gamma)~$) the set of subgroupoids (resp., set of normal sugroupoids) of $ \Gamma^{\prime},$ i.e.:\\[-0.3cm]
\begin{equation}
{\cal {S}}(\Gamma^{\prime}):= \{ H^{\prime} ~|~ H^{\prime}~ \hbox{is a {\it subgroupoid} of}~ \Gamma^{\prime}\},\label{(3.3)}\\[-0.1cm]
\end{equation}
\begin{equation}
{\cal {N}}(\Gamma^{\prime}):= \{ N^{\prime} ~|~ N^{\prime}~ \hbox{is a {\it normal subgroupoid} of}~ \Gamma^{\prime}\}.\label{(3.4)}
\end{equation}
\begin{Th}
({\bf The correspondence theorem for subgroupoids})\\
 For any surjective strong groupoid morphism $ (f , f_{0}):(\Gamma, \Gamma_{0})\longrightarrow~ (\Gamma ^{\prime}, \Gamma_{0}^{\prime})$
 there exists a one-to-one correspondence between the set $~{\cal S}(\Gamma^{\prime})~$ of the  subgroupoids of $~(\Gamma, \Gamma^{\prime}_{0})~$ and the set $~{\cal \widetilde{S}}(\Gamma)~$
 of the  sugroupoids of $~(\Gamma, \Gamma_{0})~$ which contains $~Ker(f).$
\end{Th}
\begin{Th}
({\bf The correspondence theorem for normal subgroupoids})\\
 For any surjective strong groupoid morphism $ (f , f_{0}):(\Gamma, \Gamma_{0})\longrightarrow~ (\Gamma ^{\prime}, \Gamma_{0}^{\prime})$
 there exists a one-to-one correspondence between the set $~{\cal N}(\Gamma^{\prime})~$ of the  normal subgroupoids of $~\Gamma^{\prime}~$ and the set $~{\cal \widetilde{N}}(\Gamma)~$
 of the normal sugroupoids of $~\Gamma~$ which contains $~Ker(f).$
\end{Th}
\begin{Rem}{\rm The Teorems 3.1 and 3.2 generalize the correspondence theorems for subgroups and normal subgroups by a surjective morphism of groups.}
\end{Rem}

\begin{Rem}{\rm More details and results concerning algebraic substructures of a Ehresmann groupoid,  morphisms and universal constructions  of Ehresmann groupoids  can be found by the interested reader in (\cite{giv1a, mikw88, miteza07}).}
\end{Rem}

\section{Topological groupoids and  Lie groupoids: definitions and basic results}

The research of some geometric objects associated with a topological or differential manifold sometimes leads us to a groupoid endowed with a compatible topological or differential structure, arriving at the notion of topological groupoid or Lie groupoid.

The concepts of topological groupoid and Lie groupoid are defined using the definition of the groupoid in the sense of Ehresmann.\\

{\bf 4.1. Topological groupoids}

For basic results  on topological groupoids and further references, readers should consult (\cite{brown06, ehre50, giva102, rena60}).
\begin{defn}
{\rm  A {\it topological groupoid} $~ G~$ over $~ G_{0}~$ is a groupoid $~(G,\alpha ,\beta ,\mu ,\varepsilon ,\iota, G_{0} ) $  such that $~G $  and $~G_{0} $ are topological spaces and the structure functions  are continuous, where $~G_{(2)}~$ has the subspace topology of the product topological space $ G\times G.$}
\end{defn}

If $~G~$ is a topological groupoid over $~ G_{0},~$  we will sometimes write $~(G,\alpha,\beta, G_{0}), ~(G,G_{0})~$ or $~G.~$ The map $~\alpha\times\beta: G\longrightarrow G_{0}, ~x~\to~ (\alpha\times \beta)(x) = (\alpha(x),\beta(x))~$  is called the {\it anchor map} of $~G. $
A topological groupoid $~(G, G_{0})~$ is said to be {\it transitive}, if its anchor map is surjective .

A topological groupoid $~(G,G_{0})~$ such that $~\alpha(x)=\beta(x), ~$ for all $~x\in G~$ is called {\it topological group bundle}. As example,
$~Is(G)=\{ x\in G~|~\alpha (x) = \beta(x)\}~$  is a topological group bundle, called the {\it isotropy topological group bundle} associated to $~G.$

\begin{Ex}
$(i)~~$  {\rm A topological group $~G~$ having $ e $ as unity may be considered to be a topological groupoid over $\{e\}.~$. Converesely, a topological groupoid with one
unit is a topological group.\\
$(ii)~$ Let $~{\cal PG}(G_{0})~$ be the pair groupoid} associated to $~G_{0}.~$ If $~G_{0}~$ is a topological space, then $~{\cal PG}(G_{0})~$ is a topological groupoid, called
the {\it pair topological groupoid} associated to $~G_{0}.~$
\end{Ex}
\begin{defn}
{\rm {\it (i)~} Let $~(G, \alpha, \beta,  G_{0})~$ and $~( G^{\prime}, \alpha^{\prime}, \beta^{\prime}, G_{0}^{\prime})~$ be two topological  groupoids. A morphism of
groupoids $(f, f_{0}): ((G, G_{0}) \to ( G^{\prime}, G_{0}^{\prime})~$  such that the maps $~f~$ and $~f_{0} $  are continuous, is called {\it morphism of topological groupoids} or  {\it topological groupoid morphism}.\\
{\it (ii)~} A topological groupoid morphism $~(f, f_{0}) : (G, G_{0})~\to~(G^{\prime}, G_{0}^{\prime})~$ is said to be {\it topological groupoid isomorphism }, if $~f~$ and $~f_{0}~$ are bijective maps.}
\end{defn}

If $~G_{0} = G_{0}^{\prime}~$ and $~f_{0}= Id_{G_{0}}~$ we say that $~f : G ~\to~G^{\prime}~$ is a $~G_{0}~$-{\it morphism of topological groupoids}.

The bundle of topological groupoids is defined using the notion of topological groupoid.

\begin{defn}
{\rm (\cite {miva05}) A {\it bundle of topological groupoids} is a topological fibration\\
 $\varphi=(E,p,B)~$ togheter with the structure of a topological groupoid on each fibre $~p^{-1}(b),~ b\in B.$}
\end{defn}

For the basic properties and  constructions of bundles of topological groupoids, see (\cite{miva05}).\\

{\bf 4.2. Lie groupoids}

The study of basic topics related to Lie groupoids has been addressed in a series of works, among which we mention (\cite{pradi66, mack87, giv2b, miv2a, miva06}).

We begin this subsection with the presentation of the notion of Lie groupoid and their fundamental properties.
\begin{defn}
{\rm   (\cite{mack87}). A {\it Lie groupoid $~\Gamma $ over $~M $} is a groupoid $ (G,\alpha ,\beta ,\mu ,\varepsilon ,\iota ,M) $ such that $~\Gamma $ and $~M $ are smooth manifolds,
 the maps $ \alpha $ and $~\beta $ are surjective submersions and $~\mu,~ \varepsilon,~\iota~ $ are differentiable maps.}
\end{defn}

A Lie groupoid $~(\Gamma,M $ is {\it transitive}, if the anchor  $(\alpha ,\beta ):\Gamma \to M\times M $ is a surjective submersion.

Applying properties of structure functions and common differential geometry techniques the following proposition is proved.
\begin{Prop}
(\cite{mack87}) If $(\Gamma, \alpha,\beta, \mu,\varepsilon,\iota, M)$ is a Lie groupoid over $ M, $ then:\\
$(i)~\Gamma_{{2}}~$  is a closed submanifold of the product manifold $ \Gamma\times \Gamma;$\\
$(ii) ~\mu:\Gamma_{(2)} \to \Gamma,~ (x,y) \mapsto \mu (x,y): = x y,~$ is a surjective submersion;\\
$(iii)~\varepsilon :M \to \Gamma,~$  is an injective immersion;\\
$(iv)~\iota:\Gamma \to \Gamma,$ is a diffeomorphism;\\
$(v)~ \varepsilon (M)~$  is a closed submanifold of  $~\Gamma;$\\
$(vi)~(\forall )~u \in M,$  the isotropy group $~\Gamma(u)~$ is a Lie group;\\
$(vii)~~$  If $(\Gamma, M)~$  is transitive, then  the isotropy groups $~\Gamma(u),(\forall )~u \in M,$ are isomorphic Lie groups.
\end{Prop}

A Lie group $~G $ having $~e $ as unity may be regarded  to be a  Lie groupoid over $~\{e\}.$

\begin{defn}
{\rm (\cite{mack87}) (i) Let $ (\Gamma, M)~$ and  $~(\Gamma^{\prime}, M^{\prime})~$ be two Lie  groupoids.\\
 $(i)~$ A groupoid morphism  $(f,f_{0}): (\Gamma, M)~\to ~(\Gamma^{\prime}, M^{\prime})~$  
 with the property that  $~f~$ and $~f_{0}~$ are differentiable maps, is called  {\it Lie groupoid morphism}.\\
(ii) A morphism $~(f,f_{0})~$ of Lie groupoids  such  that $~f~$ and $~f_{0}~$ are diffeomorphisms
is called  {\it isomorphism of Lie groupoids}.}
\end{defn}

A {\it Lie subgroupoid} of the Lie groupoid $~(\Gamma, M)~$ is a Lie groupoid $~(H, M^{\prime})~$  together with a Lie groupoid morphism
$~(\varphi, \varphi_{0}): (H, M^{\prime}) \to (\Gamma, M)~$  such that $~varphi: M \to \Gamma ~$ is an injective immersion.
A Lie subgroupoid  $~(H, M^{\prime})~$ is {\it wide}, if $~M^{\prime} = M~$  and $~\varphi_{0}= M. $\\

A wide Lie subgroupoid $~(H, M)~$ is a {\it normal Lie subgroupoid} of Lie groupoid $~(\Gamma, M),~$ if $~(\varphi(H), M)~$ is a normal subgroupoid
of the Ehresmann groupoid $~(\Gamma, M).$\\
Let $~(\Gamma, M)~$  be a Lie groupoid. $~Is(\Gamma)~$ is a closed submanifold of the manifold $\Gamma~$ such that canonical inclusion $~j: Is(\gamma) \to \Gamma ~$ is an injective immersion. Then $(Is(\Gamma), M)~$ is a closed normal Lie subgroupoid in $~(\Gamma, M).$ Also, if $~f: (\Gamma, M) \to (\Gamma^{\prime}, M)~$ a $~M-$ morphism of Lie groupods, then the kernel $~Ker(f)=\{x\in \Gamma~|~f(x)\in \varepsilon^{\prime}(M)\}\subseteq \Gamma~$ is a closed submanifold of the manifold $~\Gamma~$ such that the inclusion $~j:Ker(f) \to \Gamma ~$
is an injective immersion. It follows that $~(Ker(f), M)~$ is a closed normal Lie subgroupoid in $~(\Gamma, M)~$ (\cite{miteza07}).\\
A {\it Lie group bundle} is a Lie groupoid $~(\Gamma, \alpha, \beta, M)~$  such that $~\alpha (x) = \beta
(x),~(\forall )x \in \Gamma.$  For example, if $~(\Gamma, \alpha, \beta, M)~$ is a Lie groupoid, then $~Is(\Gamma)~$ is a Lie group bundle,
called the {\it isotropy group bundle} of $~\Gamma.$

\begin{Ex}
{\rm $~$  Let $~M~$ be a differentiable  manifold. The differentiable manifold $~M\times M~$ endowed with the structure of the pair groupoid $~{\cal PG}(M) = (M\times M, M)~$
is a Lie groupoid over $M,~$ called the {\it pair Lie groupoid associated to $~M~$.}}
\end{Ex}
\begin{Ex}
{\rm $~$  {\bf The direct product of two Lie groupoids.} Let $~(\Gamma, \alpha, \beta, \varepsilon, \mu, \iota, M)~$ and
 $~(\Gamma^{\prime}, \alpha^{\prime}, \beta^{\prime}, \varepsilon^{\prime}, \mu^{\prime}, \iota^{\prime}, M^{\prime})~$ be two Lie groupoids. Consider the product manifolds $~\Gamma \times \Gamma^{\prime}~$ and $~M\times M^{\prime}.~$ Then the pair $~(\Gamma \times \Gamma^{\prime}, M\times M^{\prime})~$ endowed with the groupoid structure generated by the structure functions~ $~\alpha\times \alpha^{\prime}, \beta \times \beta^{\prime}, \varepsilon \times \varepsilon^{\prime}, \mu \times \mu^{\prime}, \iota \times \iota^{\prime}~$
is a Lie gruupoid, called the {\it direct product of the Lie  groupoids} $~(\Gamma, M)$ and $~(\Gamma^{\prime}, M^{\prime}).$

The direct product of two transitive Lie groupoids is also a transitive Lie groupoid.}
\end{Ex}

A number of works are devoted to general constructions of Lie groupoids. The new Lie groupoids are obtained 
starting from one or more given Lie groupoids. Among them we list the following constructions:\\
$-~$  Lie groupoid associated with a surjective submersion and  disjoint union of a family of Lie groupoids (\cite{miv2a});\\
 $-~$ Whitney sum of two Lie groupoids over the same base, semidirect product of Lie groupoids and gauge groupoid (\cite{giv2b}); \\
$-~$ induced groupoid of a Lie groupoid (\cite{miteza07});\\
$-~$ inductive limit of an inductive system of Lie groupoids (\cite{deivan01});\\
$-~$ projective  limit of a projective  system of Lie groupoids (\cite{dege09}).

\section{The last three research directions in groupoid theory}

{\bf 5.1~ Finite groupoids}\\

A special field in groupoid theory is dedicated to the study of finite Brandt  groupoids.
 More details and results on basic properties of finite groupoids as well as their use in different areas of science, can be found in the papers (\cite{beglpt21, miv3d, dege1a, miva02, stoia2b, mapi22})

The aim of this subsection is to give some basic properties of the  finite groupoids

Let  $(G,G_0)$  be a Brandt grupoid such that  $\mid G\mid =n.$ Then the units set  $G_0$  is a finite set. If  $\mid G_0\mid =m,$  then  $1\leq m\leq n.$

A finite Brandt groupoid  $(G,G_0)$  such that  $\mid G\mid =n$ and  $\mid G_0\mid =m$ will be called a {\it groupoid of type $(n;m).$ }

If   $M$  is a finite set with  $\mid M\mid =n,$  then the pair groupoid  ${\cal{GP}}(M)$ associated to  $M$  will be denoted by  ${\cal{GP}}(n).$  For example, 
\[{\cal{GP}}(2)=\{p_{1}=(1,1),p_2=(2,2),p_3=(1,2),p_4=(2,1)\} \]  is a groupoid of type  $(4;2),$ since its units set \[ {\cal{GP}}_{0}(2)=\{p_1=(1,1),p_2
=(2,2)\}.\]  In  general ${\cal{GP}}(n)$ is a groupoid of type  $(n^{2};n).$  Also, each finite group  $G$  with $n$  elements is a groupoid of type  $(n;1).$

In the following we will refer to an important example of finite groupoid, namely the symmetric groupoid $~{\cal{S}}_n~$ , as well as to a remarkable subgroupoid of its own, namely the alternating groupoid 
$~{\cal{A}}_n~$  of  even quasipermutations  of degree  $n.$

${\bf \bullet}~~~$ {\bf The  symmetric groupoid of degree $n$}\\

Let $~({\cal S}(M),  {\cal S}_{0}(M))~$  be the groupoid of quasipermutations of the set $~M.~$ When $~M = \{ 1, 2, ..., n \},~$  we write $~{\cal S}_{n}~$ for $~{\cal S}(M)~$ and  $~{\cal S}_{n,0}~$ for $~{\cal S}_{0}(M).~$ 
 This groupoid is called  the {\it symmetric groupoid of degree} $~n.~$

The symmetric groupoid of a finite set play an important role in the study of finite groupoids, since by Cayley's theorem every finite groupoid of degree $~n~$ is
isomorphic to some subgroupoid of $~{\cal S}_{n}.~$\\

A {\it quasipermutation of length $ k~ $} of the set  $ M := \{ 1, 2, ..., n \} $ is an injective function   $f_k$ defined as:\\

\[~f_{k}=\left (\begin{array}{cccc}
i_{1} & i_{2} & \ldots & i_{k}\\
f_{k}(i_{1}) & f_{k}(i_{2}) & \ldots & f_{k}(i_{k})\\
\end{array}\right),~\]
where $~1\leq k \leq n~\quad \hbox{and}\quad ~\{ i_{1}, i_{2},
\ldots, i_{k} \}\subseteq M~$  is a ordered set  formed from $k$ elements.

\begin{Prop}
(\cite{miv3d})  Let $~n~$ be a fixed number such that $~n\geq 1.~$ The  groupoids $~{\cal S}_{n},~~{\cal S}_{n,0}~$ and $~Is({\cal S}_{n})~$ are finite groupoids of order  $~|{\cal S}_{n}|,~|{\cal S}_{n,0}|~$ resp. $~| Is({\cal S}_{n}) |,~$ where:

\begin{equation}
|~{\cal S}_{n}|~=~\sum\limits_{k=1}^{n} k! ( {n \choose k})^{2},~~
|~{\cal S}_{n,0}| ~=~2^{n} - 1 ,~~ |Is({\cal S}_{n}) | =
\sum\limits_{k=1}^{n} k! { n \choose k },\label{}
\end{equation}
(here  $~{ n \choose k}~$   denotes  the $~k-$binomial coefficient).
\end{Prop}

Let us illustrate the concepts related with the substructures of a given groupoid  in the case of the symmetric groupoid of degree $~2.~$ 

The symmetric groupoid for $ n=2~$ is  $~{\cal S}_{2} = \{f_{i}~|~ i=\overline{1,6}~\},~$ where\\

$f_{1}=\left(\begin{array}{c}

1\\

1
\end{array}\right),~~  f_{2}=\left(\begin{array}{c}

2\\

2

\end{array}\right),~~  f_{3}=\left (\begin {array}{cc}

1 & 2\\

1 & 2

\end{array}\right),~~f_{4}=\left (\begin{array}{c}

1 \\

2

\end{array}\right),$\\

$ f_{5}=\left (\begin{array}{c}

2\\

1

\end{array}\right),~~ f_{6}=\left (\begin{array}{cc}

1 & 2\\

2 & 1

\end{array}\right). $

The units set  of $~{\cal S}_{2}~$ is $~{\cal S}_{2,0}=\{ f_{1}, f_{2}, f_{3}\}.~$  Hence, $~{\cal S}_{2}~$ is a groupoid of type $~(6;3). $ The isotropy subgroupoid of $~{\cal S}_{2}~$ is $~Is({\cal S}_{2})= \{ f_{1}, f_{2}, f_{3}, f_{6}\}.~$ 

\begin{Rem}
{\rm   We have that $~| {\cal S}_{2} | = 6, ~ | Is({\cal S}_{2}) | = 4~$ and the order of $~Is({\cal S}_{2})~$ is not a divisor of $~| {\cal S}_{2} |.~$ Hence, Lagrange' s theorem for finite groups is not valid for finite groupoids.}
\end{Rem}

${\bf \bullet}~~~$ {\bf The  alternating  groupoid of degree $n$}\\

We denote the signature of quasipermutation $~f_{k}\in {\cal S}_{n} $  by $\sigma(f_{k}).$

We will say that    $~f_{k}\in {\cal S}_{n},~k\geq 2,~$  is   {\it even cvasipermutation} if  $~\sigma (f_{k})=1.~$ 

Let  $~{\cal A}_{n}~=~\{~f_{k}\in {\cal S}_{n}~|~\sigma
(f_{k})=1~ \}~$  be the set of all even quasipermutations of degree $n (n\geq 2).$

The set  $~{\cal A}_{n}~$ is a normal subgroupoid  of symmetric groupoid  $~{\cal S}_{n},~$  called the {\it alternating groupoid of degree $ n$}  (see \cite{mivan2}).

The alternating groupoid $~{\cal A}_{n}~$ contain the normal subgroupoids $~{\cal A}_{n,0}~$ and $~Is({\cal A}_{n}).~$  

\begin{Prop}
(\cite{mivan2})~ Let $~n~$ be a fixed number such that $~n\geq 2.~$ The  groupoids $~{\cal A}_{n},~~{\cal A}_{n,0}~$ and $~Is({\cal A}_{n})~$ are finite groupoids of order  $~|{\cal A}_{n}|,~|{\cal A}_{n,0}|~$ resp. $~| Is({\cal A}_{n}) |,~$ where:

\begin{equation}
|{\cal A}_{n}|~=~\displaystyle\frac{1}{2}[~\sum\limits_{k=1}^{n} k! ({ n \choose
k})^{2} - ( n^{2} - 2 n )],~~~
|{\cal A}_{n,0}| ~=~2^{n} - 1 ,\label{}
\end{equation}
\begin{equation}
|Is({\cal A}_{n})|~=~\displaystyle\frac{1}{2} [~ n +
\sum\limits_{k=1}^{n} k! { n \choose k }].\label{}
\end{equation}
\end{Prop}

\begin{Rem}
{\rm   For $ n=3,~$ we have  $~| {\cal A}_{3} | = 15  $ and   $~| {\cal A}_{3,0} | = 7.  $  Hence,  $~{\cal A}_{3} ~$ is a groupoid of type $~(15;7).$  Also  $~| Is({\cal A}_{3}) | = 9.~$ .}
\end{Rem}

In the sequel we give an example  for construction of a  finite groupoid of type $~(14;6).~$  For this, we consider the following three groupoids: the pair groupoid  ${\cal{GP}}(2), $   the symmetric groupoid $~{\cal S}_{2}~$ 
and the  group $~({\bf Z}_{4}, + ) $ of integers modulo $ 4, $ where $~{\bf Z}_{4} = \{ 0, 1, 2, 3 \}.~$

We have  $~G= \{~ a_{i}~|~i=\overline{1,14}~\},~$ where:\\
$a_{1}=p_{1},~  a_{2}=p_{2},~  a_{3}=p_{3},~  a_{4}=p_{4},$\\
$ a_{5}=f_{1}~, a_{6}=f_{2},~ a_{7}=f_{3},~ a_{8}=f_{4},~ a_{9}=f_{5},~ a_{10}=f_{6}, $\\
$a_{11}=c_{11}=0,~ a_{12}=c_{12}=1,~ a_{13}=c_{13}=2, ~a_{14}=c_{14}=3. $

Using now the groupoids  $ {\cal{GP}}(2),~{\cal S}_{2}~$ and $~{\bf  Z}_{4},~$  we consider the disjoint union of these groupoids and we denote by $~G = {\cal{GP}}(2) \coprod {\cal S}_{2} \coprod  {\bf  Z}_{4}.~$

The units set of groupoid $~G~$ is  $~G_{0} = \{ a_{1}, a_{2}, a_{5}, a_{6}, a_{7}, a_{11} \}. $\\

The multiplication  operation "$ \odot $"  defined on $~G~$  is given by:\\
$a_{i}\odot a_{j}=p_{i} .  p_{j}~$ for $ i,j\in \{1,2,3,4\},~$  if  $~p_{i}~$ and $~p_{j}~$ are composable in  $~{\cal{GP}}(2);$\\ 
$a_{k}\odot a_{l}=f_{l}\circ  f_{k}~$ for $ k,l\in \{5,6,7,8,9,10\},~$  if $~f_{k}~$ and $~f_{l}~$ are composable in $~{\cal S}_{2};$\\ 
$a_{m}\odot a_{n}=c_{m} + c_{n}~$ for all  $ m,n\in \{11,12,13,14\}~$  and $~c_{m}, c_{n}\in {\bf  Z}_{4}.~$\\ 

The structure functions  $~\alpha, \beta, \iota  $ of the groupoid $~G = \{~ a_{i}~|~i=\overline{1,14}~\}~$ are given  in the following table:

$$\begin{array}{|c|c|c|c|c|c|c|c|c|c|c|c|c|c|c|c|} \hline
a_{i}        & a_{1} & a_{2} & a_{3} & a_{4} &  a_{5} &  a_{6} & a_{7} &  a_{8} &  a_{9} &  a_{10}   &   a_{11} &   a_{12} &   a_{13} &   a_{14} \\ \hline
\alpha(a_{i}) & a_{1} & a_{2}    & a_{1} & a_{2}   &  a_{5}     &  a_{6}  & a_{7}     &  a_{5}    &  a_{6}     &  a_{7}       & a_{11}       & a_{11}       &  a_{11}  & a_{11}  \cr \hline
\beta(a_{i})  & a_{1} & a_{2}    &  a_{2}   &  a_{1}      &  a_{5}  &  a_{6}    &  a_{7}   &  a_{6}    &  a_{5}       &  a_{7}       &  a_{11}      &  a_{11}  &  a_{11}  & a_{11} \cr \hline
\iota(a_{i})      &  a_{1} & a_{2} & a_{4}  & a_{3}  & a_{5} & a_{6} & a_{7}  & a_{9}     & a_{8}  & a_{10}  & a_{11}  &  a_{14}  & a_{13} & a_{12} \cr \hline
\end{array}$$

The multiplication law  $~\odot $ on  $~G~$ is given in the following table:
{\small
$$\begin{array}{|r|c|c|c|c|c|c|c|c|c|c|c|c|c|c|} \hline
\odot & a_{1} & a_{2} & a_{3} & a_{4} & a_{5} & a_{6} & a_{7} & a_{8} & a_{9} & a_{10} & a_{11} & a_{12} & a_{13} & a_{14} \\ \hline
a_{1}  & a_{1} &  & a_{3}  &  &  &  &  &  &  &  &  &  &  &  \cr \hline
a_{2}  &  & a_{2} &  & a_{4}  &  &  &  &  &  &   &  &  &  & \cr \hline
a_{3}  &  & a_{3} &  & a_{1} &   &  &  &  &  &  &  &  &   &  \cr \hline
a_{4}  & a_{4}  &  & a_{2} &  & & & &  &   &  &  &  & &   \cr \hline
a_{5} &  &  &  &  & a_{5} & &  &a_{8}  & &  &  &  &  &  \cr \hline
a_{6}  &  &  &  &  &  & a_{6}   &  &  &  a_{9} &  & &  &  &  \cr \hline
a_{7} &  &  &  &  &  &  & a_{7} & & & a_{10} & & & &  \cr \hline
a_{8}  &  &  &  &  &  & a_{8} &  &  & a_{5}  &  &  &  &  &  \cr \hline
a_{9} &  &  &  &  & a_{9} &  & & a_{6} &  & & &  &  &  \cr \hline
a_{10}  &  &  & &  & & &  a_{10} &  &  & a_{7} & & & & \cr \hline
a_{11} &  &  &  & &  &  & & & & &  a_{11} &  a_{12} &  a_{13} &  a_{14} \cr \hline
a_{12} &  & & &  &  & & &  &  & & a_{12} & a_{13} & a_{14} & a_{11}  \cr \hline
a_{13} & &  &  &  &  & & & & &  & a_{13} & a_{14} &  a_{11} & a_{12} \cr \hline
a_{14} &  & &  &  &  &  &  & & & & a_{14} & a_{11}  &  a_{12} &  a_{13}\cr \hline
\end{array}.~$$}

The absence of the element from the arrow $~"i"~$ and the column $~"j"~$ in the table of composition law indicates the fact that the pair $~(a_{i}, a_{j})\in  G~$  is not composable. Indeed, for example we have that

$~a_{7}\odot a_{9} = f_{3}\cdot f_{5}~$ is not defined, since 
$~R(f_{3})= \{1,2\}\neq D(f_{5})=\{2\}.$\\

For to illustrate the calculation of the elements $ x \odot y $ for $ x,y \in G ~$ we give the following examples:\\

$a_{3}\odot a_{4}= p_{3}\cdot p_{4}=(1,2)\cdot (2,1)=(1,1)=p_{1}=a_{1};$\\

$a_{9}\odot a_{8} = f_{5} . f_{4}= f_{4}\circ f_{5}=\left(\begin{array}{cc}

1 \\

2

\end{array}\right )\left(\begin{array}{c}

2\\

1

\end{array}\right )=\left(\begin{array}{c}

2\\

2

\end{array}\right )= f_{2} =a_{6}.~$\\

$a_{13}\odot a_{14}= c_{13}+c_{14}= 2+3=1=c_{12}=a_{12}.$

\begin{Rem}
{\rm Let be the groupoids  $~{\cal{GP}}(m)~$  and  $~ {\bf  Z}_{n}~$  for $ m, n \geq 2.~$  Using the disjoint union and direct product of  groupoids it can be built the following  finite groupoids: \\
$~G_{1} = {\cal{GP}}(m) \coprod  {\bf  Z}_{n}    ~$ and $~G_{2} = {\cal{GP}}(m) \times  {\bf  Z}_{n}. ~$
$G_{1} $ is a groupoid of type $~(m^{2}+n; m+1)~$  and  $~G_{2} $ is a groupoid of type $~(m^{2}n; m).~$}
\end{Rem}
\begin{Rem}
{\rm As in the case of finite groups, an important theme in the theory of finite groupoids is devoted to the research of the following two topics: the problem of representing finite groupoids and the problem classifying finite-order groupoids..
  In this context, the interested reader can consult the papers (\cite{pata18, beglpt21, ibr19,  miv3d}).}\\
\end{Rem}

{\bf 5.2  Structured groupoids}

Another approach to the notion of a groupoid is that of a structured groupoid.  This  concept is obtained by adding another algebraic structure
such that the composition of the groupoid is compatible with the operation of the added algebraic structure. The most important types of structured groupoids
 are  the following: group-groupoid  (\cite{brspen, musa15, miva15}), vector groupoid (\cite{popu4a}), vector space-groupoid (\cite{miva13}),
topological group-groupoid and Lie group-groupoid (\cite{icozgu05, fana13}).

${\bf \bullet}~~~$ {\bf The  concept of group-groupoid}

The notion of group-groupoid was defined by R. Brown and Spencer in \cite{brspen}.

 In this subsection we present two important structured groupoids, namely the  concepts of  {\it group-groupoid}  and  {\it vector space-groupoid}.

In the sequel we describe  be the notion of group-groupoid as algebraic structure.  
A group structure on a nonempty set is regarded as an universal algebra determined by a binary operation, a nullary operation
 and an unary operation.

Let $ (G, \alpha, \beta, m, \varepsilon, \iota, G_{0})$ be a (Ehresmann)  groupoid. We suppose that on $G$ is defined a group structure $\omega:
G\times G \to G,~(x,y) \to \omega (x,y):=x\oplus y$. Also, we suppose that on $G_{0}$ is defined a group structure $\omega_{0}:
G_{0}\times G_{0} \to G_{0},~(u,v) \to \omega_{0}(u,v):=u\oplus v$. The unit element of $G$ (resp., $G_{0}$) is $e$ (resp.,
$e_{0}$); that is $\nu : \{\lambda\} \to G,~\lambda \to\nu(\lambda):= e $ (resp., $\nu_{0} :\{\lambda\} \to G_{0},~\lambda \to \nu_{0}(\lambda):= e_{0}$) (here $\{\lambda\}$
is a singleton). The inverse of $x\in G$ (resp., $u\in G_{0}$)  is denoted by $\bar{x}$ (resp., $\bar{u}$ ); that is $\sigma :G \to
G,~x \to \sigma (x):=\bar{x}$ (resp., $\sigma_{0}:G_{0} \to G_{0},~u \to \sigma_{0}(u):= \bar{u}$).

\begin{defn}
{\rm (\cite{brspen})  A {\it group-groupoid}  is a groupoid $ (G, G_{0})$ such that the following
conditions hold:

$(i)~~ (G, \omega, \nu, \sigma)$ and $ (G_{0}, \omega_{0}, \nu_{0}, \sigma_{0})$ are groups.

$(ii)~~$ The maps $~(\omega, \omega_{0}):(G\times G, G_{0}\times G_{0}) \to (G, G_{0}),~(\nu, \nu_{0}):(\{\lambda \}, \{\lambda \}) \to (G, G_{0}) $ and 
$ (\sigma, \sigma_{0}): (G, G_{0}) \to (G,  G_{0}) $ are groupoid morphisms.}\hfill$\Box$
\end{defn}

We shall denote a group-groupoid by $ (G, \alpha, \beta, m, \e, \iota, \oplus, G_{0})$.
\begin{Prop}
If  $~(G, \alpha, \beta, m, \varepsilon, \iota, \oplus,  G_{0})$ is  group-groupoid, then:

$(i)~ $ the multiplication $m$ and binary operation $\omega$ are compatible, that is:\\[-0.3cm]
\begin{equation}
(x\cdot y)\oplus (z\cdot t)= (x\oplus z)\cdot (y\oplus
t),~~~(\forall) (x,y), (z,t)\in G_{(2)};\label{(3.1)}
\end{equation}

$(ii)~ \alpha, \beta: (G, \oplus) \to (G_{0}, \oplus),~ i: (G,
\oplus) \to (G, \oplus)~$  and $~\varepsilon: (G_{0}, \oplus)\to (G,
\oplus) $ are morphisms of groups; i.e., for all $ x,y\in G $ and $u,v\in G_{0},$ we have:\\[-0.3cm]
\begin{equation}
\alpha(x\oplus y)= \alpha(x)\oplus \alpha(y),~~~\beta(x\oplus y)=
\beta(x)\oplus \beta(y),~~~\iota(x\oplus y)= \iota(x)\oplus \iota(y);\label{(3.2)}
\end{equation}
\begin{equation}
\e(u\oplus v)= \e(u)\oplus \e(v);\label{(3.3)}
\end{equation}

$(iii)~$ the multiplication $m$ and the unary operation $\sigma$ are compatible, that is:\\[-0.3cm]
\begin{equation}
\sigma(x\cdot y)= \sigma(x)\cdot \sigma(y),~~~(\forall) (x,y)\in
G_{(2)}.\label{(3.4)}
\end{equation}
\end{Prop}

The relation $(5.4)$ (resp., $(5.7)$) is called the {\it interchange law} between groupoid multiplication $m$ and group
operation $\omega$ (resp., $\sigma$).

We say that the group-groupoid $ (G, \alpha, \beta, m, \varepsilon, \iota,
\oplus, G_{0})$ is a {\it commutative group-groupoid}, if the
groups $G$ and $G_{0}$ are commutative.

\begin{Th}
Let $ (G, \alpha, \beta, m, \varepsilon, \iota,  G_{0})$ be a groupoid. If the following conditions are satisfied:

$(i)~~(G, \oplus)$ and $ (G_{0}, \oplus)$ are groups;

$(ii)~ \alpha, \beta: (G, \oplus)\to (G_{0}, \oplus),~\varepsilon: (G_{0},
\oplus) \to (G, \oplus)~$ and $~\iota: (G, \oplus) \to (G, \oplus) $ are morphisms of groups;

 $(iii)~$ the interchange law $(5.4)$ between  the operations $ m $ and  $ \omega $  holds,\\
 then $ (G, \alpha, \beta, m, \varepsilon, \iota, \oplus,  G_{0}) $ is a group-groupoid.
\end{Th}
According to Proposition  $5.3$ and Theorem $5.1$ we give a new definition (Definition $5.2$) for the notion of group-groupoid (this is equivalent with Definition $5.1$).
\begin{defn}
{\rm A {\it group-groupoid} is a groupoid $ (G,\alpha, \beta, m, \varepsilon, \iota, G_{0})$ such that the following conditions are satisfied:

$(i)~~ (G, \oplus)$ and $ (G_{0}, \oplus )$ are groups;

$(ii)~~ \alpha, \beta: (G, \oplus) \to (G_{0}, \oplus),~ \e:
(G_{0}, \oplus) \to (G, \oplus)~$ and $~\iota: (G, \oplus) \to (G,
\oplus) $ are morphisms of groups;

$(iii)~ $ the interchange law $(5.4)$ between  the operations $m$ and  $\oplus$ holds.} \hfill$\Box$
\end{defn}

\begin{defn}
{\rm Let $ ( G_{j}, \alpha_{j}, \beta_{j}, m_{j}, \varepsilon_{j},
\iota_{j}, \oplus_{j}, G_{j,0} ),~ j=1,2$ be two group-groupoids. A
groupoid morphism  $ (f, f_{0}): (G_{1}, G_{1,0})\to (G_{2},
G_{2,0}) $ such  that $ f $ and $ f_{0}$ are group morphisms, is
called  {\it group-groupoid morphism} or {\it morphism of
group-groupoids}.}
\end{defn}

${\bf \bullet}~~~$ {\bf The  concept of vector space-groupoid}

In paper \cite{miva13}, Mihai Ivan extends the notion of group-groupoid to the concept of {\it vector space-groupoid}.\\

Let $~(V, \alpha, \beta, m, \varepsilon, \iota, V_{0})~$ be a groupoid. we suppose that $~V~$ (resp., $ V_{0} $) is a vector space over a field $ K.$ 
For the binary operation and unary operation in the group $ V $ (resp., $ V_{0} $) we will use the notations  $~\omega:=+$ (resp., $ \omega_{0}:=+$) and $~\sigma(x):-x, x\in V$ (resp., $~\sigma_{0}(u):-=u, u\in V_{0} $).
The null vector of $ V~$ (resp., $ V_{0} $) is $ e $ (resp., $ e_{0} $ ).  The scalar multplication $~\varphi : K\times V \rightarrow V~$ (resp.., $~\varphi_{0}: K\times V_{0} \rightarrow V_{0} $) is given by:\\
$~~~~~~~~~(k,x) \mapsto \varphi(k,x):=kx~~$ (resp.,$~ (k,u) \mapsto \varphi_{0}(k, u):= k u).$\\

Consider the direct product $~(K\times V, Id\times \alpha, Id\times \beta, Id\times m, Id\times \varepsilon, Id\times \iota, K\times V_{0})~$ of the null groupoid associiated  to  $ K $  and groupoid  $ (V, V_{0} ).$
 Its set of composable elements is $~(K\times V)_{(2)} = \{ ((k_{1}, x), (k_{2},y)) \in (K\times V)^{2}~|~k_{1}=k_{2},~\beta(x)=\alpha(y)~\}.$\\
The multiplications in $~K\times V~$ is given by:\\
$~~~~~~~~~~~(k,x).(k,y):=(k, x.y),~~~~~(\forall) (x,y)\in V_{(2)},~k\in K. $\\

\begin{defn}(\cite{miva13})
{\rm A {\it vector space-groupoid} is a  groupoid $~(V, V_{0}) $ such that the following conditions hold:

$(VSG.1)~~ (V, +, \varphi )~$ and $~ ( V_{0}, +, \varphi_{0}) $ are vector spaces;\\

$(VSG.2)~~  ( V, \alpha, \beta, m, \varepsilon, \iota, +, V_{0}) $  is a commutative group-groupoid;\\

$(VSG.3)~~$ The pair $~(\varphi, \varphi_{0} ) : (K\times V, K\times V_{0} ) \rightarrow   ( V, , V_{0}) $  is a  groupoid morphism.}
\end{defn}

A vector space-groupoid is denoted by $~ ( V, \alpha, \beta, m, \varepsilon, \iota,  +, \varphi, V_{0} ) $ or $~(V, V_{0} ).$\\

In the following proposition we summarize the most important rules of algebraic calculation in a vector groupoid obtained directly from
definitions by combining the properties of structure functions of a groupoid with the properties of vector spaces.

\begin{Th}(\cite{miva13})
Let  $~( V, \alpha, \beta, m, \varepsilon, \iota, V_{0} ) $  be a groupoid. If the following conditions are satisfied:\\
$(i)~~ (V, +, \varphi ) $ and $~(V_{0}, +, \varphi_{0} ) $ are vector spaces;\\

$(ii)~~~\alpha, \beta : V \rightarrow V_{0},~\varepsilon : V_{0}\rightarrow V $ and $~\iota : V\rightarrow V ~$ are liniar maps;\\ 

$(iii)~~$ the operations $ m ~$ and $~\omega~$ are compatible, that is\\
\begin{equation}
(x\cdot y) + (z,t) = (x+z)\cdot (y+t),~~(\forall)~ (x,y), (z,t) \in V_{(2)}, \label{(3.14)}
\end{equation}
 then $~( V, \alpha, \beta, m, \varepsilon, \iota, +, \varphi, V_{0}) $  is a vector space-groupoid.
\end{Th}

The relation $~(5.8)~$ is called the {\it interchange law  between $ m ~$ and $~\omega.~$}

A fundamental results related to vector space-groupoids is contained in the characterization theorem of the vector space-groupoid concept (see Theorem 5.2).\\
\begin{Rem}
{\rm  Another extension of the algebraic concept of group-groupoid to the notion of {\bf  vector groupoid}, was defined and investigated by Vasile Popu\c ta and Gheorghe Ivan in the papers 
\cite{popu4a, popu4c}.} 
\end{Rem}
\begin{Rem} 
{\rm  The group-groupoids and their extensions(topological group-groupoids,  Lie group-groupoids, vector space-groupoid) are mathematical structures that have proved
to be useful in many areas of science, see for instance \cite{akiz18, guari17, guric18, popu4b, mivde18, mann15}.}
\end{Rem}

{\bf  5.3. Applications of  groupoids}\\

The study of Ehresman groupoids, topological groupoids, Lie groupoids and Lie algebroids is motivated by their applications in various branches of mathematics, namely in:\\
- algebra and  category theory(\cite{avmp20, baga22,  deis93, dugiv94, higgi71, mapi22};\\
-  topology and algebraic topology (\cite{ brown06, dugiv94, giva1c92, west71});\\
-analysis, locally compact groupoids and quantum mechanics (\cite{goste06, giva102, puta5a, miteza07, rare01, west71});\\
- theory of representations of topological groupoids (\cite{ hedrea5a, mann15, rena60, wein96});\\
- differential geometry and geometric mechanics (\cite{puta5b, miteza07, mack87, pope07, pradi66});\\
- symplectic geomeetry and Hamiltonian mechanics (\cite{ cdw87, migopr09, mikw88, wein87});\\
-combinatorics (\cite{ziva06});\\
- crystallography (\cite{john00}) ;\\
- study of some dynamical systems in Hamiltonian form using Lie groups, Lie groupoids, Leibniz algebroids and Lie algebroids (\cite{mana06, comgimg09, coma04, ivop06, opris5a, migopr10});\\
- as well as in the theory of classical and fractional dynamical systems on Lie groupoids and Lie algebroids (\cite{gmio06, imod07, miva23, migi18, migo09}).

Moreover, various types of constructions and classes of actions of Lie groupoids on differential manifolds are used in the construction of principal bundles with structural Lie groupoids (\cite{giv5b, giv5c}). These principal bundles constitute an important 
study tool in various applications of differential geometry in theoretical mechanics (\cite{miv32,  miteza07}). For example, connections on principal bundles with structural Lie groupoids are used to prove integrability conditions of Lie algebras or
 Lie algebroids. The linear connections on Lie algebras or Lie algebroids  have many applications in Hamiltonian mechanics and other areas of physics  (\cite{dege5a, opris5b, marti01, sahi14}).\\

\vspace*{0.2cm}

Author's adress\\[0.2cm]
\hspace*{0.7cm}West University of Timi\c soara. Seminarul de Geometrie \c si Topologie\\
\hspace*{0.7cm} Department of Mathematics\\
\hspace*{0.7cm} Bd. V. P{\^a}rvan,no.4, 300223, Timi\c soara, Romania\\
\hspace*{0.7cm}E-mail: gheorghe.ivan@e-uvt.ro\\


{\small
\begin{thebibliography}{99}

\bibitem{akiz18} H.F. Akiz, A note on  the fundamental  groupoid as a vector-space groupoid. {\it Academic Researches in Math. and Sci.}, Ankara, First Edition,  2018,  21-26.

\bibitem{mana06} M. Anastasiei, Mechanical systems on Lie algebroids. {\it Algebras Groups Geometry}, 23(3) (2006),  235-245.

\bibitem{avmp20} J. Ávila, V. Marín and  H. Pinedo, Isomorphism theorems for groupoids
and some applications,{\it Int. J. Math. Math. Sci.}, {\bf 2020}, no. 3967368,
1–10, 2020. Doi:10.1155/2020/3967368

\bibitem{baga22}  M.A. Babayo, G.U. Garbab,  Alternative  proof of the isomorphism theorem of groups.
{\it Journal of Equations Algebras}, {\bf 10}(1)(2022),  72-75. 

\bibitem{beglpt21}  G. Beier, C. Garcia, W. G. Lautenschlaeger, J. Pedrotti,  T. Tamusiunas, Generalizations of Lagrange and Sylow theorems for groupoids. {\it Arxiv:}2101.07420v1 [math.RA], 19 Jan  2021, 1-19.

\bibitem{brand} H. Brandt, {\it \"{U}ber eine Verallgemeinerung der Gruppen-Begriffes}. Math. Ann., {\bf 96} (1926), 360-366.

\bibitem{brown} R. Brown, \textit{From Groups to Groupoids: a brief survey.} Bull. London Math. Soc., \textbf{19} (1987), 113-134.
https//doi.org/10.1112/blms/19.2.113.

\bibitem{brspen} R. Brown, C.B. Spencer, {\it G-groupoids, crossed modules and the fundamental groupoid of a topological group}. Proc. Kon. Nederl.
Akad. Wet., {\bf 79} (1976), 296-302.

\bibitem{brown06} R. Brown, {\it Topology and groupoids}. BookSurge LLC, U.K., 2006.

\bibitem{comgimg09} D. Comănescu, M. Ivan, Gh. Ivan, Lyapunov stability for the e-revised dynamics of the rigid body with three linear controls. Journal of Computational Analysis and Applications, 11(3) (2009), 546-559. www.pdf.academia.edu.

\bibitem{cdw87}  A. Coste, P. Dazord, A. Weinstein, {\it Groupoïd$\acute{e}$s symplectiques}. Publ. Dept. Math. Lyon, 2/A
(1987),1-62.

\bibitem{coma04} J. Cortes, E. Martinez, Mechanical control systems on Lie algebroids. {\it IMA J.
Math. Control Inform.}, {\bf 21} (4)(2004), 457–492.

\bibitem{deis93} M. Deaconescu, Gh. Ivan, G. Silberberg,  Proceedings of the International Conference on Group Theory: Timi\c soara, September, 1992. {\it An. Univ. Timi\c soara}, Ser. Științ.. Mat., 1993, Special Issue, 1–192.

\bibitem{dege5a} M. Degeratu, M. Ivan, Linear connections on Lie algebroids. {\it Proceed. of the 5th  Conf. of Balkan Society of Geometers} (2005-Mangalia). Balkan Press, Bucharest, 2006, 44-53. www.Linear-connections-on-Lie-algebroids/pdf.

\bibitem{dege1a} M. Degeratu, G. Ivan, M. Ivan, On the cyclic subgroupoids of a Brandt groupoid. {\it Proceed. Int. Conf. Comp., Commun. Control}
(ICCCC 2006-Oradea). Vol. I, Supplementary Issue (2006), 181-186. www. journal.univagora.ro; doi.org/10.1587/ijccc.2006.5.

\bibitem{deivan01} M. Degeratu,  Gh. Ivan, {\it Inductive limits of Lie groupoids}.  An. Univ. Timi\c soara Ser. Mat.-Inform. 39 (2001), Mathematics 175-187. pdf.researchgate/publication/254036373.

\bibitem{dege09} M. Degeratu,  Projective limits of Lie groupoids. An. Univ. Oradea, Fasc. Mat. 16 (2009), 183–189.

\bibitem{dugiv94} B. Dumons, Gh. Ivan, Introduction à la théorie des groupoïdes. Dept. Math. Univ. Poitiers (France),URA,C.N.R.S. D1322,  86, 1994, 1- 86.

\bibitem{ehre50} C. Ehresmann, \textit{O\'{e}uvres Compl\'{e}tes. Parties I.1, I.2. Topologie Alg\'{e}brique et G\'{e}ometrie Diff\'{e}rentielle}.
Dunod, Paris, 1950.

\bibitem{fana13} M.R. Farhangdoost, T. Nasirzade, {\it Geometrical categories of generalized Lie groups and Lie group-groupoids}. Iran J Sci Technol A., (2013), 69-73.

\bibitem{goste06} M. Golubitsky, I. Stewart, {\it Nonlinear dynamics of networks: the groupoid formalism}. Bull. Amer. Math.
Soc., {\bf 43}(3)(2006),  305-364.

\bibitem{guari17} M.H. G$\ddot{u}$rsoy, H. Aslan, I. $\dot{I}$cen, {\it Generalized crossed modules and group-groupoids}. Turk. J. Math. {\bf 41}(6) (2017),  1535-1551.

\bibitem{guric18} M.H. G$\ddot{u}$rsoy, I. $\dot{I}$cen, {\it Coverings of structured Lie groupoids}. Hacettepe J. Math. Statistics {\bf 47}(1) (2018), 845-854.


\bibitem{hedrea5a}. C. Hedrea, I. Gole\c t, R. Negrea, Geometric prequantization of the dual a Lie algebra. {\it Int.J. Geom. Meth. Modern Phys.}, {\bf 10}(6)(2013), 1320009.

\bibitem{higgi71} P. J. Higgins, Notes on Categories and Groupoids, {\it Math. Studies} {\bf 32}, London, 1971.

\bibitem{ibr19}  A. Ibort, M. Rodriguez, {\it On the structure of finite groupoids and their representations}.  Symmetry, {\bf 11}(3) (2019), 4-14.
https://doi.org/10.3390/sym110304.

\bibitem{icozgu05} I. $\dot{I}$cen, A. F. $\ddot{O}$zcan,  M. H. G$\ddot{u}$rsoy, \emph{Topological group-groupoids and their coverings}.
Indian J. Pure Appl. Math., {\bf 36}(9)(2005), 493-502.

\bibitem{john00}  C.K. Johnson, Crystallographic groups, groupoids  and orbifolds. {\it  Workshop  on Orbifolds, Groupoids and Their Applications}, University of Wales, Bangor, September, 2000.

\bibitem{giv1a} Gh. Ivan, Strong morphisms of groupoids. {\it Balkan J.Geom. Appl.}(BJGA), {\bf 4} (1)(1999), 91-102.
 https://www.emis.de/journals/BJGA/v04n1/B04-1-IVAN.pdf.

\bibitem{giv5b} Gh. Ivan, Principal fibre bundles with structural Lie groupoid. {\it Balkan J. Geom. Appl.}( BJGA), {\bf 6}(2) (2001), 39-48.

\bibitem{giv1b} Gh. Ivan, Algebraic constructions of Brandt groupoids. {\it Proceed. of the Algebra Symposium}, " Bab\c s- Bolyai" University, Cluj-Napoca, (2002), 69-90. pdf.px-pict.com, www.07-ivan-pdf.

\bibitem{giv2b} Gh. Ivan, Extensions of Lie Groupoids. {\it Recent  Advances In Geometry and Topology}.  Cluj Univ. Press, 2004, 199-210.

\bibitem{giv5c} Gh. Ivan, Actions of Lie groupoids on fibred manifolds. {\it Tensor (N.S.)}, {\bf 66}(2) (2005), 142-156.

\bibitem{giva102} Gh. Ivan,   Haar systems of measures on locally compact groupoids. Proceed. Conf. Nat. de Analiză şi Aplicaţii, Timişoara-2000.  An. Univ. Timişoara (2002), 74-90.

\bibitem{giva1c92} Gh. Ivan,  On cohomology of Brandt groupoids. Proceed. of Int. Conf. on Group Theory , Timișoara-1992, (1992), 17-20.

\bibitem{ivop06} Gh. Ivan, D. Opri\c s, Dynamical systems on Leibniz algebroids. {\it Differ. Geom.-Dyn. Syst.}, {\bf 8}(2006), 127-137.

\bibitem{popu4b} Gh. Ivan,  V. Popu\c ta, Vector subgroupoids. The computation of all vector subgroupoids of a modular vector groupoid. {\it Carpathian J. Math.}, {\bf 28}(2) (2012), 265–270.

\bibitem{puta5a} Gh. Ivan, M. Puta, Groupoids and some of their applications in geometry and quantum mechanics. {\it Proceed. of the Int. Conf. on Group Theory} (Timi\c soara,1992). An. Univ. Timi\c soara Ser. \c Stiin\c t. Mat. 1993, Special Issue, 40-41.

\bibitem{puta5b} Gh. Ivan, M. Puta, Groupoids, Lie$-$Poisson structures and quantization. {\it Proceed. of the Int. Conf. on Group Theory} (Timi\c{s}oara, 1992), An. Univ. Timi\c soara Ser. \c Siin\c t. Mat. 1993, Special Issue, 117–123.

\bibitem{opris5a} Gh. Ivan, M. Ivan, D. Opri\c s, Structures on affine Lie algebroids and applications. {\it Contemporary Geometry and Topology and Related Topics}, 173–185, Cluj Univ. Press, Cluj-Napoca, 2008. https://www.imar.ro/~purice/Inst/Applications-cluj.php.

\bibitem{gmio06} Gh. Ivan, M. Ivan, D. Opri\c s, Fractional Euler-Lagrange and fractional Wong equations for Lie algebroids. {\it Proceedings of the 4-th International Colloquium Mathematics and Numerical Physics},  Bucharest-2006,  73-80.

\bibitem{imod07} Gh. Ivan, M. Ivan, D. Opri\c s, Fractional dynamical systems on fractional Leibniz algebroids, {\it An. \c Stiin. Univ., "Al. I. Cuza" Ia\c si} (S.N.), Mat., {\bf 53} (2007), Supl., 222-234.

\bibitem{miva02} M. Ivan,  On the symmetric groupoid. Seminar Arghiriade 32, Timişoara, Univ. de Vest din Timişoara, Facultatea de Matematică, 2002, 1-23.

\bibitem{mivan2}  M. Ivan, {\it The alternating groupoid $~{\mathcal A}_{n}.~$} Proceedings of the $10^{th}$ International Symposium of Mathematics and its Applications,
" Politehnica " University of Timi\c soara-2003, (2003), 227-234. https://www.mat.upt. ro/Simpozion Matematica/index.html

\bibitem{miv3d} M. Ivan, General properties of the symmetric groupoid of a finite set. {\it An. Univ. Craiova, Ser. Mat.Comp. Sci. Ser.}, {\bf 30}(2) (2003),  109–119.
https://aucmcs.ro/index.php.

\bibitem{miv2a} M. Ivan, Bundles of Lie groupoids. {\it Recent  Advances In Geometry and Topology}. Cluj Univ. Press, 2004,  211–219.
www.Bundles-of-Lie-groupoids.pdf

\bibitem{miv32} M. Ivan, Refinements of a principal bundle with Lie groupoid structure. Tensor N.S., 66(2) ( 2005 ), 157–166.

\bibitem{miva05} M. Ivan,  Bundles of topological groupoids, {\it Universitatea din Bac\u au. Studii \c si Cercet\u ari Ştiin\c tifice, Seria Matematic\u a},
nr.{\bf 15} (2005),  43–54. www.ub.ro/pubs/scssm.

\bibitem{miva06} M. Ivan, Lie's Functor for bundles of Lie groupoids, {\it Tensor (N.S.)}, {\bf 67}(3) (2006),  294–304. tensor-ns@nifty.ne.jp

\bibitem{miteza07} M. Ivan, {\it Contributions to the study of groupoids} (Contribu\c tii la studiul grupoizilor). PhD thesis. West University of Timi\c soara. Faculty of Mat.-Inform., 2007, 1-126.

\bibitem{miva13}  M. Ivan,  Vector space-groupoids, {\it Eur. J. Pure Appl Math.}, {\bf 6}(4) (2013), 469-484.
www.eyup.+1a-Ivan.pdf.

\bibitem{miva15}  M. Ivan, Algebraic properties of ${\cal G}$-groupoids. {\it ArXiv:}1512.09012v1
[math.GR] 30  Dec 2015, 1-11. www.1512.09012.pdf

\bibitem{miva23} M. Ivan, Analysis stability for fractional Volterra model with two controls. {\it International Journal of Applied Mathematics}, {\bf 36}(6) (2023), 801-814. Doi.org/10,12732/ijam.v36i6.5.

\bibitem{mivde18} M. Ivan, M. Degeratu, A method for computing of group-subgoupoids of finite group-groupoids. {\it An. Univ. Oradea}, Fasc. Matematic\u a, {\bf 25}(1) (2018),  25-32.
www.researchgate.net/publication.

\bibitem{migi18}  M. Ivan,  Gh. Ivan, On the fractional Euler top system with two parameters.  {\it Int. J. Modern Eng. Research}, {\bf 8}(4)(2018), 10-22.

\bibitem{opris5b} M. Ivan,  D. Opri\c s, Hamiltonian formalism on the transformation Lie algebroid. {\it An. Univ. Vest Timiş. Ser. Mat.-Inform.} {\bf 46}(1) (2008),  85–95. www.Hamiltonian-formalism-on-the-transformation-Lie-algebroid.pdf.

\bibitem{stoia2b}.  M. Ivan, G. Stoianov, Determination of the transitive subgroupoids of a groupoid. An. St. Univ. de Vest, Timişoara, Seria Mat.- Inform, {\bf 43}(2) (2005), 73–87.

\bibitem{migo09}  M. Ivan,  Gh. Ivan, D.  Opri\c s, Fractional equations of the rigid body on the pseudo-orthogonal group $ SO(2,1) $. {\it Int. J. Geom. Met. Mod. Phys.}, {\bf 6}(7)(2009), 1181-1192.  Doi:10.1142/S0219887809004168.

\bibitem{migopr09} M. Ivan, Gh. Ivan, D. Opriş, Numerical integration for the equations of the heavy top in  Hamiltonian form. {\it Tensor (N.S.)}, {\bf 71}(1) (2009), 51-50. 

\bibitem{migopr10} M. Ivan, Gh. Ivan, D. Opriş, The Euler-Arnold equations in Hamiltonian form of heavy top. {\it Carpathian J. Math.}, {\bf 26}(2) (2010), 202-209. 

\bibitem{mack87} K. Mackenzie, \textit{Lie Groupoids and Lie Algebroids in Differential Geometry.} London Math. Soc.,Lecture Notes Series,
\textbf{124}, Cambridge Univ.Press., 1987.

\bibitem{mann15} M.K. Mann'a, Groupoids in topological vector space, JP J. Geom. Topol., 17 (2015), 127-144.

\bibitem{mapi22} V. Marín  and H. Pinedo, Groupoids: direct products, semidirect products and solvability, {\it Algebra and Discrete Math.}, {\bf 33}(2) (2022), 92-107. Doi:10.12958/adm 1772.

\bibitem{marti01}  E. Martinez, Geometric formulation of mechanics on Lie algebroids. In: {\it Proceed. of the VIII-th Fall Workshop on Geometry
and Physics} (Medina del Campo, 1999). Publicaciones de la RSME, \textbf{2} (2001), 209-222.

\bibitem{mikw88}  K. Mikami and A. Weinstein, Moments and reduction for symplectic groupoid actions, {\it Publ. RIMS Kyoto Univ.}, {\bf 42} (1988), 121-140.

\bibitem{musa15} O. Mucuk O, T. Sahan and N. Alemdar, Normality and quotients in crossed modules and group-groupoids, {\it Appl. Categor. Struct.}, {\bf 23} (2015),415-428.

\bibitem{pata18}  A. Paques, T. Tamusiunas,  The Galois correspondence theorem for groupoid actions. {\it J. of Algebra}, {\bf 509} (2018), 105-123.
https://doi.org/10.1016/j.jalgebra.2018.04.034.

\bibitem{pope07} L.  Popescu, Aspects of Lie algebroid geometry and hamiltonian formalism. {\it An. \c Stiin\c t. Univ. Al. I. Cuza Ia\c si}, Mat. (N.S.), {\bf 53} (2007), Supl., 297-308.

\bibitem{popu4a} V. Popu\c ta, Gh. Ivan, Vector groupoids. {\it  Theoretical Mathematics and Applications} (TMA), International Scientific Press, {\bf 2}(2) (2012), 1-12.

\bibitem{popu4c} V. Popu\c ta, Gh. Ivan, Morphisms of vector groupoids. {\it Communications in Math. and  Applications}, {\bf 4}(1) (2013),  61–73. http://www.rgnpublications.com.

\bibitem{pradi66} J. Pradines, {Th}\'{e}orie de Lie pour les groupo\"{\i}des diff\'{e}rentiables. {\it  Comptes Rendus Acad. Sci. Paris},
\textbf{263 }(1966), 907-910.

\bibitem{rare01} A. Ramsey, J. Renault, {\it Groupoids in Analysis, Geometry and Physics}. Contemporary Mathematics, {\bf 282}, AMS Providence, RI, 2001.

\bibitem{rena60} J. Renault, \textit{A groupoid approach to C*-algebras.} Lecture
Notes Series, \textbf{793}, Springer-Verlag, 1980.

\bibitem{sahi14}F. \c Sahin,  Lie groupoids and generalized contact manifolds. {\it Abstract and Applied Analysis},{\bf 2014} (2014),Article ID,   270-715.
https//doi.org/10.1155/2014/270715.

\bibitem{wein87} A. Weinstein, Symplectic groupoids and Poisson manifolds. {\it Bull. Amer. Math. Soc.}, {\bf 16} (1987), 101-104.

\bibitem{wein96} A. Weinstein, Groupoids: Unifying internal and external symmetries. {\it Notices Amer. Math. Soc.}, \textbf{43} (1996), 744-752.

\bibitem{west71} J. Westman, {\it Groupoid Theory in Algebra, Topology and Analysis}. University of California at Irvine, 1971.

\bibitem{ziva06} R. T. \v{Z}ivaljević, \textit{Groupoids in combinatorics - applications of a theory of local symmetries.} Preprint,
arXiv:math.CO / 0605508, 2006, 1-23.

\end{thebibliography}}
\end{document}